\documentclass[11pt,twoside]{report} 
\usepackage{amsmath,amsthm,amssymb,amsfonts}
\usepackage{graphicx,epsf,epsfig}
\usepackage{setspace}
\theoremstyle{plain}
\newtheorem{thm}{Theorem}[section]

\newtheorem{cor}[thm]{Corollary}
\newtheorem{lem}[thm]{Lemma}

\theoremstyle{definition}
\newtheorem{defn}[thm]{Definition}
\newtheorem{exe}[thm]{Example}

\theoremstyle{remark}

\newtheorem*{note}{Note}

\numberwithin{equation}{section}

\newcommand{\Tri}{{\bigtriangleup}}
\newcommand{\To}{{\rightarrow}}

\newcommand{\Z}{{\mathbb N}}

\newcommand{\Q}{{\mathbb Q}}
\newcommand{\R}{{\mathbb R}}
\newcommand{\PPR}{{\mathbb P^2 \mathbb R}}
\newcommand{\RRT}{\widetilde{\R^2}}
\newcommand{\fsum}{\hat{+}}
\newcommand{\bsum}{\check{+}}
\newcommand{\mylim}[2]{{\displaystyle{\lim_{#1 \To #2}}~}}
\newcommand{\mylimsup}[2]{{\displaystyle{\limsup_{#1 \To #2}}~}}
\newcommand{\myliminf}[2]{{\displaystyle{\liminf_{#1 \To #2}}~}}
\newcommand{\emptypage}{\clearpage{\pagestyle{empty}\cleardoublepage}}

\begin{document}

\title{Two-Dimensional Analogs of the Minkowski $?(x)$ Function}
\author{Andrew Marder\\[.25in]
        Department of Mathematics\\
        Williams College\\
        Williamstown, MA 01267\\
        amarder@wso.williams.edu\\[.25in]}
\date{\today}

\pagenumbering{roman}               
\maketitle                          
\emptypage                          

\pagestyle{plain}                   

\tableofcontents                    
\listoffigures                      
\emptypage                          

\pagenumbering{arabic}                                 

\chapter{Background}
One of the key properties of continued fractions is that a real number $\alpha$ is 
quadratic irrational if and only if it has an eventually periodic continued fraction 
expansion \cite{Khinchin}.  This property motivated Hermann Minkowski to define his 
question-mark function
\[
? : [0,1] \To [0,1],
\]
which has the following characteristics:
\begin{itemize}
\item $?(x)$ is strictly increasing and continuous.
\item If $x=p/q$ is rational, then $?(x)=r/2^s$ is a pure dyadic number.
\item If $x$ is a quadratic irrational, then $?(x)$ is rational.
\item The inverse image of the rational numbers is exactly the set of quadratic
irrationals.
\end{itemize}
In addition to being continuous and monotonically increasing, $?(x)$ also has 
derivative zero almost everywhere, making $?(x)$ a naturally occurring example of a 
singular function \cite{Girgensohn, Derivative, NewLight, Salem}.  Our aim 
is to extend Minkowski's function to two dimensions.  Ideally, we would like to find
functions from $\R^2$ to $\R^2$ that are strictly increasing, continuous, map cubic 
irrationals to rationals\footnote
{
  There is no known multidimensional continued fraction with the property that a real
  number $\alpha$ is cubic irrational if and only if its multidimensional continued
  fraction expansion is eventually periodic.  Therefore we cannot be certain of
  sending all cubic irrational numbers to rationals.  Instead, we will be satisfied 
  with mapping a subset of cubic irrationals to a subset of rationals.
}, 
and have derivative zero almost everywhere.

Following the model of the question-mark function, we will construct two similar
functions, both stemming from work in \cite{fareybary}.  Our functions will map a 
two-dimensional simplex 
(a triangle) to itself.  The first map will be defined by partitioning the triangle, 
initially using a ``weighted Farey'' partition, and then by a ``weighted Bary'' 
partitioning.  The function $\delta$ will map the weighted Farey 
triangle to the weighted Bary triangle.  We will see that the infinite sequence of 
Farey subtriangles containing a given point can be used as a multidimensional 
continued fraction.  Specifically, we show that points with eventually periodic Farey 
sequences are at worst cubic irrational.  And 
we show that periodicity in the Bary partition implies rationality.  Therefore our 
weighted Farey-Bary map will map a natural class of cubic irrational points to a 
natural class of rationals.  To finish our discussion of the weighted Farey-Bary map, 
we prove an analog of singularness by showing that, almost everywhere, the area of the
Bary triangles approaches zero far more quickly than the area of the Farey triangles. 

The third chapter develops a new version of the Farey-Bary map.  Like the original, 
this function will carry a subset of cubic irrational points to a subset of rational 
points.  Unlike its predecessor, this revised Farey-Bary map is continuous.  This and 
other interesting properties suggest that this new function is the natural extension 
of the question-mark function.  The final section of this paper mimics Salem's proof
that $?'(x)=0$ a.e. to prove an analog of singularness for the continuous Farey-Bary 
map.

I would like to thank Thomas Garrity, my advisor from Williams College, and Olga 
Beaver, my second reader.  Without their work and guidance, this senior thesis could
never have happened.

\section{A Review of Minkowski's Question-Mark Function}

The Minkowski question-mark function is defined inductively.  We set the initial 
values
\[
?(0) = 0 \text{ and } ?(1) = 1.
\]
Now suppose we know the values of $?(p_1/q_1)$ and $?(p_2/q_2)$.  We then set
\[
?( \frac{p_1+p_2}{q_1+q_2} ) = 
         \frac{ ?(\frac{p_1}{q_1}) + ?(\frac{p_2}{q_2}) }{2}.
\]
Using this definition, we know the value of $?(x)$ for any rational number $x$, and by
continuity arguments we can determine the value of $?(x)$ for any real number 
$x \in [0,1]$.  Although this is an elegant definition it is rather difficult to
work with.  The following is an alternative definition whose framework we will follow
closely when defining our Farey-Bary maps.  We now define two partitions of the unit
interval.  The first, the Farey partition, will be the domain of the question-mark 
function, and the second, the Bary partition, will be the range.

\begin{defn}
The \textit{Farey sum} of two rational numbers $p_1/q_1$ and $p_2/q_2$, each in lowest
terms, is defined to be
\[
\frac{p_1}{q_1} \fsum \frac{p_2}{q_2} = \frac{p_1+p_2}{q_1+q_2}.
\]
\end{defn}

\begin{defn}
The \textit{Bary sum} of two rational numbers $x$ and $y$ is defined to be the average
of the two numbers
\[
x \bsum y = \frac{x+y}{2}.
\]
\end{defn}

We produce two sequences of partitions, $I_k$ and $\tilde{I}_k$, of the unit interval.
Each partition will split the unit interval into $2^k$ subintervals, and the $k^{th}$
partition will be a refinement of the previous $(k-1)^{st}$ partition.  Both sequences
begin with the unit interval
\[
I_0 = \tilde{I}_0 = [0,1].
\]
Given the partition $I_{k-1}$, suppose that the endpoints of each of its $2^{k-1}$
subintervals are rational numbers.  We form the next partition $I_k$ by taking the 
Farey sum of the endpoints of $I_{k-1}$.  The first few partitions in this sequence are
\begin{eqnarray*}
& I_0 = \left[ \frac{0}{1}, \frac{1}{1} \right] & \\
& I_1 = \left[ \frac{0}{1}, \frac{1}{2}, \frac{1}{1} \right] & \\
& I_2 = \left[ \frac{0}{1}, \frac{1}{3}, \frac{1}{2}, 
               \frac{2}{3}, \frac{1}{1} \right] & \\
& I_3 = \left[ \frac{0}{1}, \frac{1}{4}, \frac{1}{3}, \frac{2}{5}, \frac{1}{2}, 
               \frac{3}{5}, \frac{2}{3}, \frac{3}{4}, \frac{1}{1} \right] & 
\end{eqnarray*}

The sequence of partitions $\tilde{I}_k$ is constructed in the same manner as the 
$I_k$, but we replace the Farey sum with the Bary sum.  
\begin{eqnarray*}
& \tilde{I}_0 = \left[ \frac{0}{1}, \frac{1}{1} \right] & \\
& \tilde{I}_1 = \left[ \frac{0}{1}, \frac{1}{2}, \frac{1}{1} \right] & \\
& \tilde{I}_2 = \left[ \frac{0}{1}, \frac{1}{4}, \frac{1}{2}, 
               \frac{3}{4}, \frac{1}{1} \right] & \\
& \tilde{I}_3 = \left[ \frac{0}{1}, \frac{1}{8}, \frac{1}{4}, \frac{3}{8},\frac{1}{2},
               \frac{5}{8}, \frac{3}{4}, \frac{7}{8}, \frac{1}{1} \right] &
\end{eqnarray*}

\begin{figure}[hbtp]
\centering
\includegraphics[width=\textwidth]{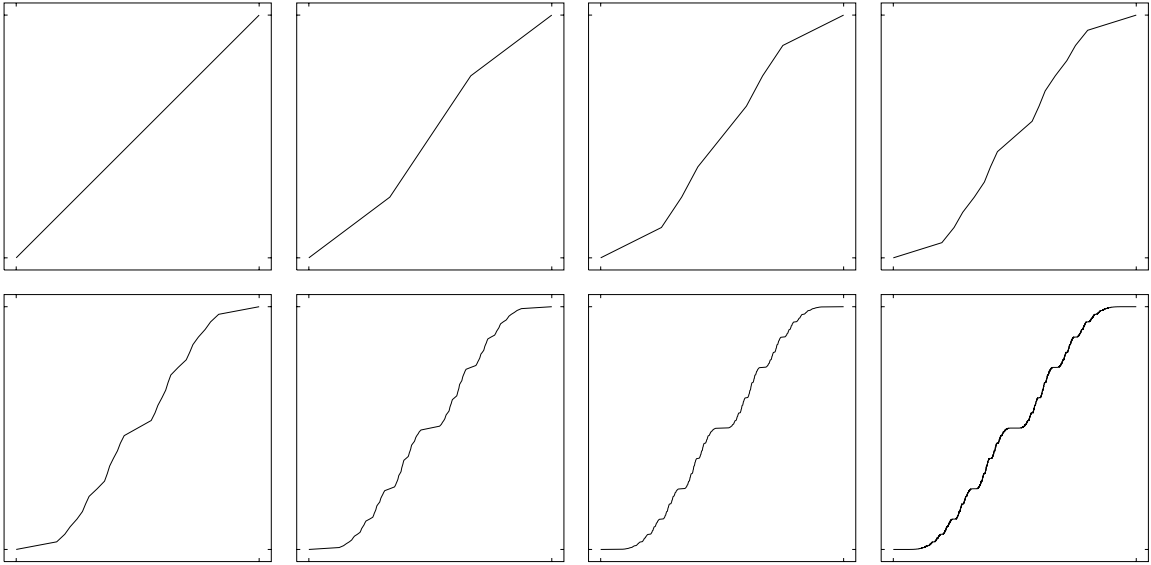}
\caption{The progression of the $?_k(x)$.}
\label{Minkowski}
\end{figure}

We define a sequence of functions $\{ ?_k(x) \}$, where each $?_k(x)$ maps the 
endpoints of $I_k$ to the corresponding endpoints of $\tilde{I}_k$ and is extended 
linearly to the rest of the
unit interval.  Some of these functions can be seen in Figure \ref{Minkowski}.  It 
turns out that this sequence $\{ ?_k(x) \}$ converges uniformly to the 
Minkowski question-mark function $?(x)$.

As stated earlier, $?(x)$ is singular, meaning it is continuous, strictly increasing,
and has derivative zero almost everywhere.  Using our partitions, $I_k$ and 
$\tilde{I}_k$, we can recast this derivative condition in terms of lengths of 
subintervals.  Fix $\alpha \in [0,1]$.  For each $k$, let $i_k \in I_k$ be the 
subinterval that contains $\alpha$, and let $\tilde{i}_k \in \tilde{I}_k$ be the
corresponding subinterval in the Bary partition.  As proven by Salem in \cite{Salem},

\begin{thm}
For almost all $\alpha \in [0,1]$ if the following limit exists and is finite then
\[
\mylim{k}{\infty} \frac{length(\tilde{i}_k)}{length(i_k)} = 0.
\]
\end{thm}

This theorem provides the most natural language for us to generalize the 
failure of differentiability to higher dimensions.  For our functions 
$\delta : \R^2 \To \R^2$ we will show a similar condition using the areas of 
corresponding subtriangles.

\section{The Domain}\label{RRT}
There are two reasons to phrase Farey and Bary sums in terms of matrices.  
First, using matrices simplifies our proofs concerning periodic Farey and Bary 
sequences.  Second, it will be relatively easy to calculate the area of a 
triangle given its matrix representation, which will be important in proving 
singularness.  In order to make use of matrix multiplication we will work in the 
following space
\[
\RRT = \R \times \R \times \R^* / \{(x,y,z) \sim (\lambda x, \lambda y, \lambda z), 
       \forall \lambda \in \R^* \},
\]
where $\R^* = \R - \{0\}$.  $\RRT$ is a subset of $\PPR$ that we will identify with
$\R^2$.

An element $v \in \RRT$ can be visualized as a line in $\R^3$ that passes 
through the origin but does not lie in the $XY$ plane.  We will define the $Class(v)$ 
to be all points on the line $v$, excluding the origin.  If $(dx, dy, 1)$ is the slope
of the line $v$ then 
\[
Class(v) = \{ \lambda (dx,dy,1) : \lambda \in \R^* \}.
\]
We will restrict our focus to lines that have rational slope because they will be
equivalent to points in $\R^2$ that have rational coordinates.

\begin{defn}
Define $\phi : \R \times \R \times \R^* \To \R^2$ to be
\[
\phi (x,y,z) = (\frac{x}{z}, \frac{y}{z}).
\]
\end{defn}

Given an element $v \in \RRT$, we want to show that $\phi$ sends each point in 
$Class(v)$ to the same point in $\R^2$.  Let $(dx,dy,1)$ be the slope of $v$ and 
$\lambda \in \R^*$.  Then
\[
\phi( \lambda \cdot dx, \lambda \cdot dy, \lambda ) = 
( \frac{\lambda \cdot dx}{\lambda}, \frac{\lambda \cdot dy}{\lambda} ) = 
(dx,dy) = \phi(dx,dy,1).
\]
Hence $\phi$ is also a function from $\RRT$ to $\R^2$, giving rise to the subsequent
revision.

\begin{defn}
Let $v \in \RRT$ and $(x,y,z) \in Class(v)$.  We redefine 
\[
\phi : \RRT \To \R^2
\]
by
\[
\phi(v) = (\frac{x}{z},\frac{y}{z})
\]
\end{defn}

\begin{thm}
$\phi$ is a bijection.
\end{thm}

\begin{proof}
First, we show that $\phi$ is injective.  Let $u,v \in \RRT$ with 
$\phi(u) = \phi(v)$, we want to show $u = v$.  Let $(dx,dy,1)$ be the slope of $u$ and
 $(ds,dt,1)$ be the slope of $v$ then
\[
\phi(u) = (dx,dy) = (ds,dt) = \phi(v),
\]
implying the two lines have the same slope and therefore must be equal.

Next, we have to show that $\phi$ is surjective.  Given $(x,y) \in \R^2$, 
there exists a line $v \in \RRT$ that passes through the origin and the point 
$(x,y,1) \in \R \times \R \times \R^*$, thus $\phi(v) = (x,y)$.
\end{proof}

For the rest of this paper, we will use $\phi$ implicitly.  When 
dealing with points in the Farey partition we will use $(p,q,r)$ to represent a point
in $\R^2$, where 
$p,q,r \in \Z$.  For points in the Bary partition we will use $(x,y,1)$, with 
$x,y \in \Q$.  $\phi$ gives us the freedom to use these representations 
interchangeably.

We now define a metric $\tilde{d}$ on $\RRT$.  Given two elements $u,v \in \RRT$ we
define $\tilde{d}$ by
\[
\tilde{d}(u,v) = d(\phi(u),\phi(v)),
\]
where $d$ is the Euclidean metric on $\R^2$.  Notice that 
\begin{eqnarray*}
\tilde{d}( (\lambda x,\lambda y, \lambda),(x,y,1) ) 
     &=& d( \phi(\lambda x,\lambda y, \lambda),\phi(x,y,1) )\\
     &=& d( (x,y), (x,y) )\\
     &=& 0.
\end{eqnarray*}  
We will use this metric while working with limits in $\RRT$ to prove a well-known 
result concerning eigenvectors.

\section{Farey and Bary Sums}

When constructing the Minkowski question-mark function we defined the Farey and Bary 
sums of two rational numbers.  Similarly, we must define the Farey and Bary sums of 
two rational points in the plane.  

\begin{defn}
Let $v_1, v_2 \in \R^2$, and
\begin{eqnarray*}
v_1 &=& (\frac{p_1}{r_1}, \frac{q_1}{r_1}),\\
v_2 &=& (\frac{p_2}{r_2}, \frac{q_2}{r_2}),
\end{eqnarray*}
where, for each $i$, the $p_i,q_i$, and $r_i$ share no common factor.  The 
\textit{Farey sum} $\hat{v}$ of the $v_i$ is then
\[
\hat{v} = v_1 \fsum v_2 = ( \frac{p_1+p_2}{r_1+r_2}, \frac{q_1+q_2}{r_1+r_2} ).
\]
\end{defn}

One bonus of this definition is that it scales up nicely
\begin{eqnarray*}
v_1 \fsum v_2 \fsum v_3 &=& (v_1 \fsum v_2) \fsum v_3 \\
                &=& (\frac{p_1+p_2}{r_1+r_2},\frac{q_1+q_2}{r_1+r_2}) \fsum 
(\frac{p_3}{r_3},\frac{q_3}{r_3}) \\
                &=& (\frac{p_1+p_2+p_3}{r_1+r_2+r_3},\frac{q_1+q_2+q_3}{r_1+r_2+r_3}).
\end{eqnarray*}

To define the Farey sum in $\RRT$, let $v_1,v_2 \in \RRT$ and
\begin{eqnarray*}
v_1 &=& (p_1,q_1,r_1),\\
v_2 &=& (p_1,q_1,r_1),
\end{eqnarray*}
where, for each $i$, the $p_i,q_i$, and $r_i$ share no common factor.  The Farey sum 
is a straightforward vector sum
\[
\hat{v} = v_1 \fsum v_2 = (p_1,q_1,r_1) + (p_2,q_2,r_2) = (p_1+p_2,q_1+q_2,r_1+r_2).
\]
Verifying this definition, we see
$\phi(v_1 \fsum v_2) = (\frac{p_1+p_2}{r_1+r_2},\frac{q_1+q_2}{r_1+r_2})
= \phi(v_1) \fsum \phi(v_2)$.

Our Farey-Bary map wouldn't be complete with just a Farey sum.  We now define the Bary
sum.

\begin{defn}
Let $v_1,...,v_n \in \R^2$.  The \textit{Bary sum} $\check{v}$ of the $v_i$ is an 
average of the $n$ points
\[
\check{v} = v_1 \bsum ... \bsum v_n = \frac{v_1 + ... + v_n}{n} =
(\frac{x_1+...+x_n}{n},\frac{y_1+...+y_n}{n}).
\]
\end{defn}

To place this definition in $\RRT$, let $v_1,...v_n \in \RRT$, and
\begin{eqnarray*}
v_1 &=& (x_1,y_1,1), \\
 &\vdots& \\
v_n &=& (x_n,y_n,1).
\end{eqnarray*}
The Bary sum $\check{v}$ of the $v_i$ is an average of the $n$ vectors
\[
\check{v} = v_1 \bsum ... \bsum v_n
          = \frac{v_1 + ... + v_n}{n}
          = ( \frac{x_1 + ... + x_n}{n}, \frac{y_1 + ... + y_n}{n}, 1).
\]

We have finished all the necessary background material and are prepared to develop the
weighted Farey-Bary map.

\chapter{The Weighted Farey-Bary Map}
In this chapter, we generalize the Farey-Bary map as defined in \cite{fareybary} by 
associating weights with each vertex of
the triangle.  We find that the weighted Farey partition fails to fully partition the 
triangle if the weight on the third vertex is greater than one.  But, if the weight on
the third vertex is equal to one then previous results of the Farey-Bary map hold 
(a class of cubic irrationals is mapped to a class of rationals and the function is 
singular).  This is a surprising asymmetry for a function that, at first glance, 
appears symmetric.

\section{Weighted Farey and Bary Sums}

Before we can define our function we need to incorporate weights into the Farey and 
Bary sums.  For both definitions we want the result of summing $n$ rational points to 
be a rational point.  

\begin{defn}
Given $m_1,...m_n \in \Z^+$ that share no common factor, and $v_1,...,v_n \in \RRT$, 
with
\begin{eqnarray*}
v_1 &=& (p_1, q_1, r_1),\\
 &\vdots& \\
v_n &=& (p_n, q_n, r_n),
\end{eqnarray*}
we define the \textit{weighted Farey sum} $\hat{v}$ of the $v_i$ to be
\begin{eqnarray*}
\hat{v} &=& m_1 v_1 \fsum ... \fsum m_n  v_n \\
        &=& (m_1p_1 + ... + m_np_n, m_1q_1 + ... + m_nq_n, m_1r_1 + ... + m_nr_n).
\end{eqnarray*}
\end{defn}

Note that we have adjusted this definition from the original definition of Farey sum.
We no longer require that the $p_i,q_i,$ and $r_i$ share no common factor.  This 
subtlety will be touched upon later but is of little consequence.

\begin{defn}
Given weightings $w_1,...,w_n \in \Q^+$ satisfying $w_1+...+w_n = 1$, and 
$v_1,...,v_n \in \RRT$ with
\begin{eqnarray*}
v_1 &=& (x_1,y_1,1),\\
 &\vdots& \\
v_n &=& (x_n,y_n,1),
\end{eqnarray*}
we define the \textit{weighted Bary sum} to be
\begin{eqnarray*}
\check{v} &=& w_1 v_1 \bsum ... \bsum w_n v_n \\
          &=& (w_1x_1 + ... + w_nx_n, w_1y_1 + ... + w_ny_n, 1).
\end{eqnarray*}
\end{defn}

\section{The Weighted Farey and Bary Partitions}

\begin{figure}[p]
\centering
\includegraphics[height=.9\textheight]{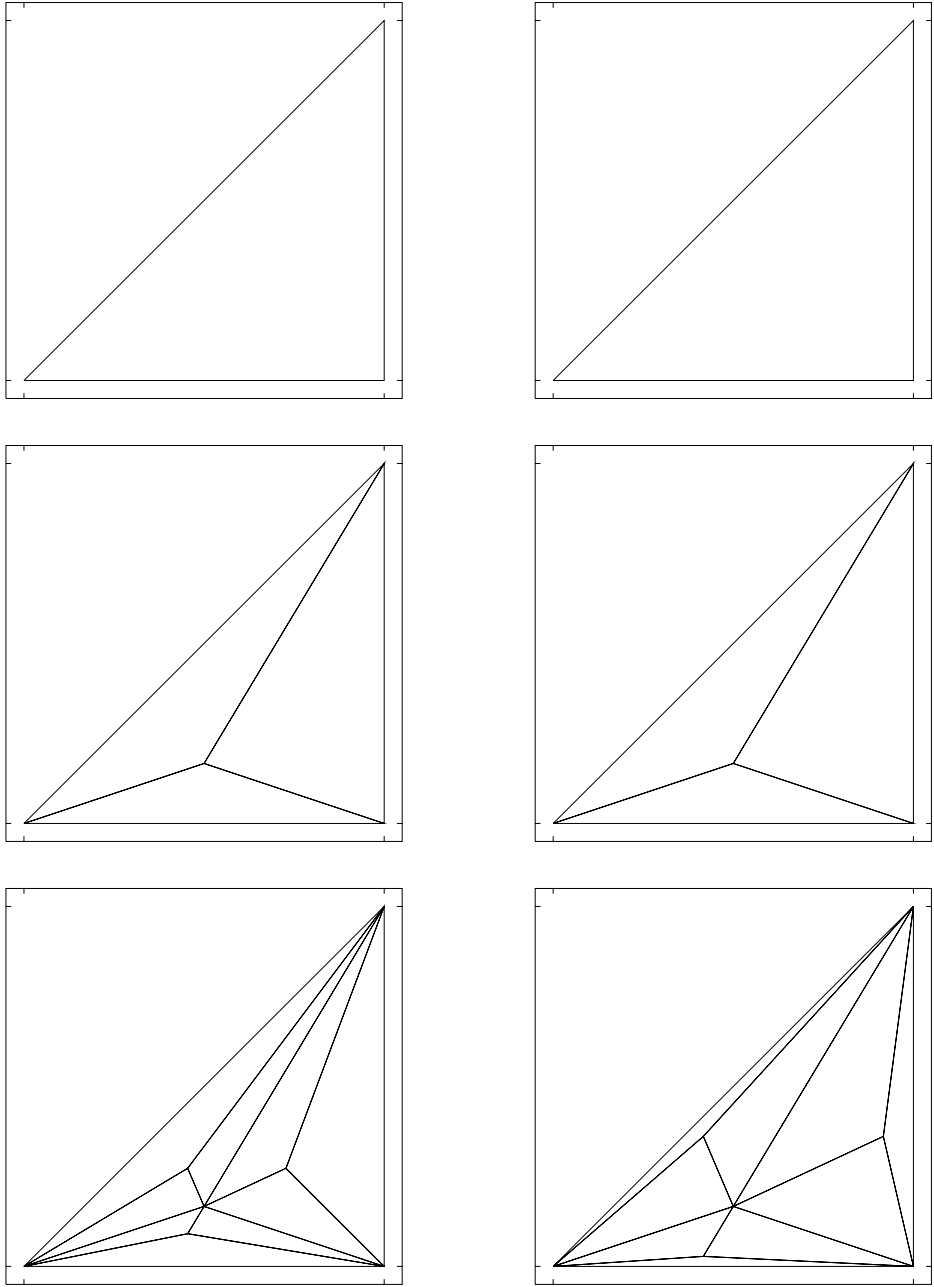}
\caption[The progression of the weighted Farey and Bary partitions.]
        {The weighed Farey and Bary partitions with weights $m_1=3,m_2=2,m_3=1$.}

\label{weightedpartitions}
\end{figure}

In this section we define two partitions of the triangle
\[
\Tri = \{ (x,y) : 1 \ge x \ge y \ge 0 \}.
\]
The first partition will yield the domain of our desired function, while the second 
will yield the range.  Figure \ref{weightedpartitions} shows the first six stages of 
the weighted Farey and Bary partitions using weights $m_1 = 3, m_2 = 2,$ and $m_3=1$.
As some of these technical details can be difficult to visualize this diagram will
clarify much of this section.

Given natural numbers $m_1,m_2,m_3$ that share no common factor, we will define the 
weighted Farey partition of the triangle by a sequence of partitions $\{ P_n \}$ such
that each $P_n$ will consist of $3^n$ subtriangles of $\Tri$ and each $P_n$ will be
a refinement of the previous $P_{n-1}$.  Let $P_0$ be the initial triangle $\Tri$.  
The three vertices of $\Tri$ are $(0,0),(1,0),(1,1) \in \R^2$.  Since we are working 
in $\RRT$, we have
\begin{eqnarray*}
v_1 &=& (0,0,1), \\
v_2 &=& (1,0,1), \\
v_3 &=& (1,1,1). 
\end{eqnarray*}
Taking the weighted Farey sum of these vertices, produces a point on the interior of
$P_0$
\begin{eqnarray*}
\hat{v} &=& m_1 v_1 \fsum m_2 v_2 \fsum m_3 v_3 \\
        &=& (m_2 + m_3,m_3,m_1 + m_2 + m_3) \\
        &=& (\frac{m_2 +m_3}{m_1+m_2+m_3}, \frac{m_3}{m_1+m_2+m_3}).
\end{eqnarray*}
This point refines $P_0$ by cutting the triangle into three new subtriangles.  Thus 
determining the partition $P_1$.  

We define the rest of the $P_n$ inductively.  Suppose we have the partition $P_n$ that determines $3^n$
subtriangles of $\Tri$.  We split each of these triangles into three smaller 
subtriangles as follows.  Given a subtriangle in $P_n$ that has vertices 
$v_1,v_2,v_3$, we compute the weighted Farey sum of the three vertices producing an 
interior point $\hat{v}$, which divides the triangle into three new subtriangles.  
Calculating the weighted Farey-center for each subtriangle of $P_n$, gives us the next
partition $P_{n+1}$ of $\Tri$.  We denote the full partitioning of $\Tri$ by $\Tri_F$
and call it the weighted Farey partitioning.

\begin{note}
In both the weighted Farey and Bary partitions it is important to number the vertices
of a subtriangle consistently because each vertex is assigned a different weight.  In 
the weighted Farey partition, we will order the vertices such that
\[
r_1 \leq r_2 \leq r_3,
\]
where $v_i = (p_i,q_i,r_i)$.  Although this is not a major concern now, it will be an
important fact in later proofs.
\end{note}

We define another sequence of partitions $\{ \tilde{P}_n \}$ of $\Tri$ such that each 
$\tilde{P}_n$ will consist of $3^n$ subtriangles of $\Tri$ and each $\tilde{P}_n$ is
a refinement of the previous $\tilde{P}_{n-1}$.  The only difference here is that we 
will use the weighted Bary sum.  To compute this sum we need rational weights that sum
to one.  To match the proportions of the weighted Farey partition we will use
\begin{eqnarray*}
w_1 &=& \frac{m_1}{m_1+m_2+m_3},\\
w_2 &=& \frac{m_2}{m_1+m_2+m_3},\\
w_3 &=& \frac{m_3}{m_1+m_2+m_3}.
\end{eqnarray*}
We call this partitioning of $\Tri$ the weighted Bary partitioning and denote it 
$\Tri_B$.

\begin{figure}[p]
\centering
\includegraphics[height=.9\textheight]{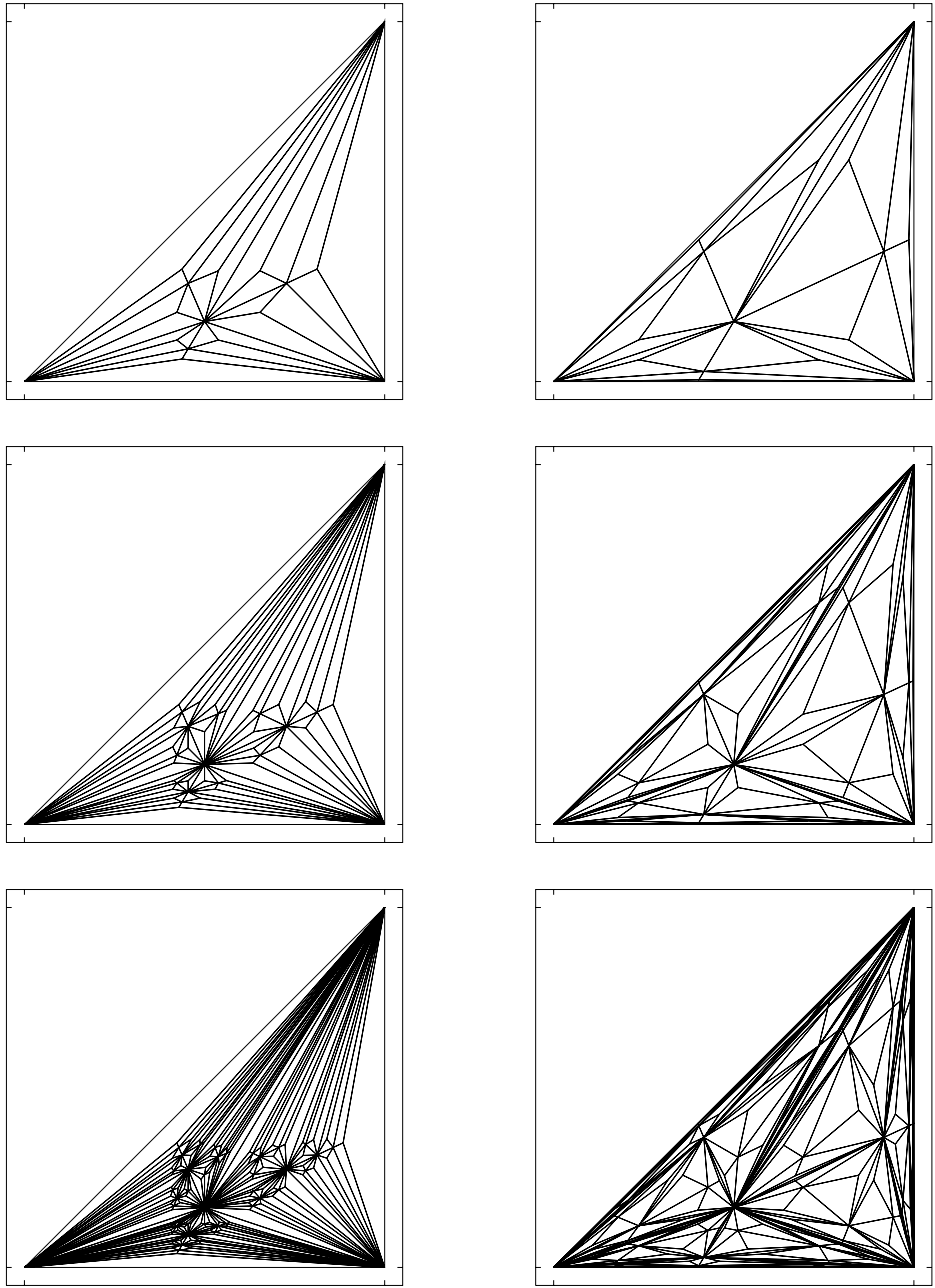}
\caption[The weighted Farey and Bary partitions continued.]
        {The weighted Farey partitions (left) and 
	 the weighted Bary partitions (right).}
\end{figure}

\section{The Weighted Farey-Bary Map}

We now define a sequence of functions $\{ \delta_n \}$ that will be the basis of the
weighted Farey-Bay map.  Keep in mind that when $m_1=m_2=m_3=1$ the weighted 
Farey-Bary map is identical to the original map in \cite{fareybary}.  First, we 
introduce some new notation.  Our two partitions 
$P_n$ and $\tilde{P}_n$ determine subtriangles of $\Tri_F$ and $\Tri_B$.  The 
expression $\langle v_1(n), v_2(n), v_3(n) \rangle$ will denote a general subtriangle 
of $P_n$ with vertices $v_1(n), v_2(n),$ and $v_3(n)$.  To denote a general
subtriangle of $\tilde{P}_n$ we will use 
$\langle \tilde{v}_1(n), \tilde{v}_2(n), \tilde{v}_3(n) \rangle$.

\begin{defn}
Given a natural number $n$, define $\delta_n$ to send any vertex in the $n^{th}$ Farey
partition $P_n$ to the corresponding vertex in the $n^{th}$ Bary partition
$\tilde{P}_n$.  That is, define $\delta_n$ on any subtriangle 
$\langle v_1(n), v_2(n), v_3(n) \rangle$ of $P_n$ by
\[
\delta_n(v_i(n)) = \tilde{v}_i(n).
\]
Finally, for any point $(x,y)$ in the subtriangle with vertices 
$\langle v_1(n), v_2(n), v_3(n) \rangle$, set
\[
\delta_n(x,y) = \alpha \tilde{v}_1(n) + \beta \tilde{v}_2(n) + \gamma \tilde{v}_3(n),
\]
where
\[
(x,y) = \alpha v_1(n) + \beta v_2(n) + \gamma v_3(n).
\]
\end{defn}
Since the point $(x,y)$ is in the interior of the subtriangle we have that
\[
\alpha + \beta + \gamma = 1
\]
and
\[
0 \leq \alpha, \beta, \gamma < 1.
\]

We define part of the weighted Farey-Bary map as follows.  If, for some $n$, $v$ falls
on an edge of a partition triangle in $P_n$ then we set $\delta(v) = \delta_n(v)$.  We
are mapping the edges of each Farey subtriangle to the corresponding edges in the Bary
partition.  Treatment of points that do not lie on any edge will be dealt with in 
section \ref{definition}.

\section{Farey Iteration as Multidimensional Continued Fraction}

Minkowski's question-mark function maps quadratic irrational numbers to rationals.  
Like the original Farey-Bary map, we
hope to send a class of cubic irrational points to a class of rational points.  
Keeping this target in mind we begin defining the Farey sequence.

There are two types of points in the Farey partition.  Those that land on an
edge of a partition subtriangle and those that do not.  We will focus on the latter.

\begin{defn}
A point $v \in \Tri$ is an \textit{interior point} of the Farey partition 
if it does not land on an edge of any partition subtriangle.
\end{defn}

Since the number of subtriangles is countable we see 
that almost all points in $\Tri$ are interior points.  The Farey sequence of a point
$v$ will be a numeric representation of the nested sequence of Farey subtriangles 
that contain $v$.  If $v$ is interior then it will have a unique Farey sequence.

To define the Farey sequence formally we first relate the vertices of the 
$(n-1)^{st}$ subtriangle that contains $(\alpha, \beta)$ with the vertices of the 
subtriangle at the next stage.  Suppose that 
\[
(\alpha, \beta) \in \langle v_1(n-1), v_2(n-1), v_3(n-1) \rangle 
                    \subseteq P_{n-1}.
\]
Applying the next partition, we decompose the triangle into three new 
subtriangles.  If we let $\langle v_1(n), v_2(n), v_3(n) \rangle$ denote the 
subtriangle into which $(\alpha, \beta)$ falls, we see that there are three 
possibilities for the vertices of $\langle v_1(n), v_2(n), v_3(n) \rangle$

\begin{figure}[hbtp]
\centering
\includegraphics[height=2.5in]{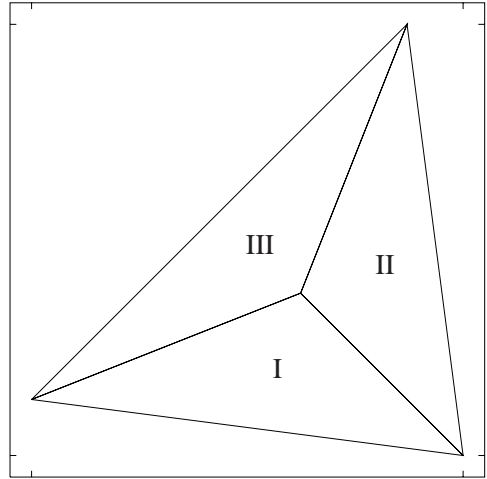}
\caption{The three Farey subtriangles.}
\end{figure}

\begin{enumerate}
\centering
\item $\langle v_1(n-1), v_2(n-1), m_1 v_1(n-1) \fsum m_2 v_2(n-1) 
       \fsum m_3 v_3(n-1 \rangle$
\item $\langle v_2(n-1), v_3(n-1), m_1 v_1(n-1) \fsum m_2 v_2(n-1) 
       \fsum m_3 v_3(n-1 \rangle$
\item $\langle v_1(n-1), v_3(n-1), m_1 v_1(n-1) \fsum m_2 v_2(n-1) 
       \fsum m_3 v_3(n-1 \rangle$
\end{enumerate}

For each $(\alpha, \beta) \in \Tri$, we now associate a sequence of positive integers.
This sequence of integers will be determined by the sequence of Farey subtriangles 
that contain $(\alpha, \beta)$.  

Although our eventual notation may seem unmotivated, it will be central in our proof 
of singularness.  Consider the following three possibilities.  Start with a triangle
whose vertices, $v_1,v_2,$ and $v_3$, maintain the convention that 
$r_1 \le r_2 \le r_3$.  Suppose we perform $L$ type I operations.  The new triangle 
will have vertices
\[
v_1, v_2, (1+...+m_3^{L-2})(m_1v_1 \fsum m_2v_2) \fsum m_3^{L-1}(m_1v_1 
           \fsum m_2v_2 \fsum m_3v_3). \footnote
{
We are abusing notation here.  When $L=1$ we say $(1+...+m_3^{L-2})=0$.
}
\]
If we perform a type II operation on the triangle, and then $L-1$ type I 
operations, the new triangle will have vertices
\[
v_2, v_3, (1+...+m_3^{L-2})(m_1v_2 \fsum m_2v_3) \fsum m_3^{L-1}(m_1v_1 
           \fsum m_2v_2 \fsum m_3v_3).
\]
And if we perform a type III operation, and then $L-1$ type I 
operations, the new triangle will have vertices
\[
v_1, v_3, (1+...+m_3^{L-2})(m_1v_1 \fsum m_2v_3) \fsum m_3^{L-1}(m_1v_1 
           \fsum m_2v_2 \fsum m_3v_3).
\]
We have set the framework for the following notation.

Define the Farey sequence $\{ a_1(i_1), a_2(i_2), ... \}$ to be such that each $a_k$ 
is a positive integer and each $i_k$ represents either case I,II, or III.  The value 
of $a_k(i_k)$ denotes the action of first applying a type $i_k$ operation and then 
$a_k - 1$ type I operations.

Note that by the time we have the $k^{th}$ term of the sequence we
have performed $s_k = a_1 + ... + a_k$ Farey operations on 
$\Tri$.  We associate to each $(\alpha, \beta) \in \Tri$ the sequence that yields the 
corresponding sequence of Farey subtriangles that contain $(\alpha, \beta)$.  If 
$(\alpha,\beta)$ is an interior point then its Farey sequence will be unique.

Now we place our three Farey transformations in terms of matrix multiplication.  Again
we will use $\langle v_1(n), v_2(n), v_3(n) \rangle$ to denote a subtriangle of 
$P_n$.  Our initial triangle is written as
\[
M_0 =
( v_1(0) ~ v_2(0) ~ v_3(0) )
=
\left(
\begin{array}{ccc}
0 & 1 & 1 \\
0 & 0 & 1 \\
1 & 1 & 1 \\
\end{array}
\right) .
\]
And we will write subsequent triangles as
\[
M_n = 
( v_1(n) ~ v_2(n) ~ v_3(n) )
=
\left(
\begin{array}{ccc}
p_1(n) & p_2(n) & p_3(n) \\
q_1(n) & q_2(n) & q_3(n) \\
r_1(n) & r_2(n) & r_3(n) \\
\end{array}
\right).
\]
Our three transformations become:
\begin{enumerate}
\centering
\item 
$M_n = M_{n-1} 
\left(
\begin{array}{ccc}
1 & 0 & m_1 \\
0 & 1 & m_2 \\
0 & 0 & m_3 \\
\end{array}
\right)$

\item 
$M_n = M_{n-1} 
\left(
\begin{array}{ccc}
0 & 0 & m_1 \\
1 & 0 & m_2 \\
0 & 1 & m_3 \\
\end{array}
\right)$

\item 
$M_n = M_{n-1} 
\left(
\begin{array}{ccc}
1 & 0 & m_1 \\
0 & 0 & m_2 \\
0 & 1 & m_3 \\
\end{array}
\right)$
\end{enumerate}

\begin{thm}
All possible $M_n$ have determinant $\pm 1$ if and only if $m_1 = m_2 = m_3 = 1$.
\end{thm}

\begin{proof}
This follows immediately from observing that $det(M_0) = 1$ and that the determinants
of the transformations matrices are also $\pm 1$ if and only if $m_1 = m_2 = m_3 = 1$.
\end{proof}

\begin{thm}
If $m_1 = m_2 = m_3 = 1$ then the entries in each column vector of $M_n$ share no 
common factor.
\end{thm}

\begin{proof}
From the previous theorem we know that $M_n$ must have determinant $\pm 1$.  Without 
loss of generality, assume that the entries of the third column share a common factor.
Then
\[
M_n
=
\left(
\begin{array}{ccc}
p_1(n) & p_2(n) & p_3(n) \\
q_1(n) & q_2(n) & q_3(n) \\
r_1(n) & r_2(n) & r_3(n) \\
\end{array}
\right)
=
\left(
\begin{array}{ccc}
p_1(n) & p_2(n) & c \cdot p'_3(n) \\
q_1(n) & q_2(n) & c \cdot q'_3(n) \\
r_1(n) & r_2(n) & c \cdot r'_3(n) \\
\end{array}
\right),
\]
where all the values are positive integers and $c > 1$.  Taking the determinant we 
find
\[
det(M_n) 
= 
c \cdot 
\left|
\begin{array}{ccc}
p_1(n) & p_2(n) & p'_3(n) \\
q_1(n) & q_2(n) & q'_3(n) \\
r_1(n) & r_2(n) & r'_3(n) \\
\end{array}
\right|, 
\]
which does not equal $\pm 1$ because it is a product of two integers one of which 
is greater than one.  Concluding, the entries of each column share no common factor.
\end{proof}

This is a useful result for our unweighted function, but it points out a 
slight issue when using more interesting weights.

\begin{exe}
Here is an example of a non-trivial weighting scheme that results in a column whose 
entries share a common factor of two.  Let $m_1 = 3, m_2 = 2,$ and $m_3 = 1$, consider
the following sequence of transformations
\[
\left(
\begin{array}{ccc}
0 & 1 & 1 \\
0 & 0 & 1 \\
1 & 1 & 1 \\
\end{array}
\right)
\left(
\begin{array}{ccc}
1 & 0 & 3 \\
0 & 0 & 2 \\
0 & 1 & 1 \\
\end{array}
\right)
\left(
\begin{array}{ccc}
0 & 0 & 3 \\
1 & 0 & 2 \\
0 & 1 & 1 \\
\end{array}
\right)
\left(
\begin{array}{ccc}
1 & 0 & 3 \\
0 & 1 & 2 \\
0 & 0 & 1 \\
\end{array}
\right)
=
\left(
\begin{array}{ccc}
1 & 3 & 14 \\
1 & 1 & 8  \\
1 & 6 & 26 \\
\end{array}
\right).
\]
The factor of two places extra weight on the third vertex when calculating the 
weighted Farey sum.  Because the result of a weighed Farey sum may have a common
factor we had to adjust our definition of the weighted Farey sum to accommodate. 
However this adjustment poses no problem in our analysis.
\end{exe}

\subsection{Farey Periodicity Implies Cubic Irrationality}

We are poised to prove the major result of this section.  First, we 
prove an important lemma.

\begin{lem}\label{eigenvector}
Let $A$ be a $3 \times 3$ real matrix.  If one of the columns $C_i(k)$ of $A^k$
converges to a vector $v \in \RRT$ then $v$ is an eigenvector of $A$.
\end{lem}

\begin{proof}
Using the metric defined in section \ref{RRT} we take the limit of the sequence
$\{ A C_i(k) \}$
\begin{eqnarray*}
\mylim{k}{\infty} A C_i(k) &=& Av  \\
\mylim{k}{\infty} C_i(k+1) &=& Av  \\
                         v &=& Av. \\
\end{eqnarray*}
Since this equality occurs in $\RRT$ there must exist a 
$\lambda \in \R^*$ such that $Av = \lambda v$.
\end{proof}

The following theorem is the first half of showing that the weighted Farey-Bary map 
sends a natural class of cubic points to a natural class of rational points.

\begin{thm}\label{cubic}
Suppose that $\{ a_k(i_k) \}$ is an eventually periodic Farey sequence that converges 
to the point $(\alpha,\beta)$.  Then $\alpha$ and $\beta$ are algebraic numbers with 
$deg(\alpha) \le 3$, $deg(\beta) \le 3$ and
\[
dim_{\Q}\Q[\alpha,\beta] \le 3.
\]
\end{thm}

\begin{proof}
There are two facts central to this proof.  First, as seen in \cite{TriSequences}, if
$(a ~ b ~ 1)^T$ is an eigenvector of a $3 \times 3$ matrix with rational entries
then
\[
dim_{\Q}\Q[a,b] \le 3.
\]
The second being the preceding lemma.

Recall that the vertices of the partition subtriangles containing $(\alpha,\beta)$ 
correspond to the columns of their matrix representations, which are products of the 
various transformation matrices.  Having assumed $\{ a_k(i_k) \}$ is eventually 
periodic, denote the product of the initial non-periodic matrices $A$ and the product 
of the periodic part $B$.  Then some of the Farey partition triangles about the point 
$(\alpha, \beta)$ are given by
\[
A,AB,AB^2,AB^3,...
\]
By assumption, the columns of the matrices $A,AB,AB^2,...$ converge to a multiple
of $(\alpha ~ \beta ~ 1)^T$.  Thus the columns of the matrices $B,B^2,B^3,...$ must 
converge to a multiple of $A^{-1}(\alpha ~ \beta ~ 1)^T$.  Applying the previous lemma
we know that $A^{-1}(\alpha ~ \beta ~ 1)^T$ must be an eigenvector of $B$.  
Therefore $\alpha$ and $\beta$ must have the desired properties.
\end{proof}

\section{Iteration in the Barycentric Range}

To finish proving that $\delta$ sends a subset of cubic irrational points 
to a subset of rational points, we will prove that 
points with eventually periodic Bary sequences must be rational.  First, we define the
weighted Bary sequence $\{ \tilde{a}_k(i_k) \}$ as we did the Farey sequence.  The 
value of $\tilde{a}_k(i_k)$ denotes the action of first applying a type $i_k$ Bary 
transformation and then $\tilde{a}_k - 1$ type I operations.  Again we translate these
operations in terms of matrix multiplication.

\begin{enumerate}
\centering

\item 
$\tilde{M}_n = \tilde{M}_{n-1} 
\left(
\begin{array}{ccc}
1 & 0 & w_1 \\
0 & 1 & w_2 \\
0 & 0 & w_3 \\
\end{array}
\right)$

\item 
$\tilde{M}_n = \tilde{M}_{n-1} 
\left(
\begin{array}{ccc}
0 & 0 & w_1 \\
1 & 0 & w_2 \\
0 & 1 & w_3 \\
\end{array}
\right)$

\item 
$\tilde{M}_n = \tilde{M}_{n-1} 
\left(
\begin{array}{ccc}
1 & 0 & w_1 \\
0 & 0 & w_2 \\
0 & 1 & w_3 \\
\end{array}
\right)$
\end{enumerate}

Recall that $w_i = m_i/(m_1+m_2+m_3)$ for $i=1,2,$ or $3$.  Therefore these three 
matrices are stochastic, meaning that their columns sum to one.  This fact will
be critical in proving the following proof.

\subsection{Bary Periodicity Implies Rationality}

From lemma \ref{eigenvector} we know that if the columns of $A^k$ converge to a vector
$v$ then $v$ must be an eigenvector of $A$.  To prove that Bary periodicity implies
rationality we must first show that if $A$ is a stochastic matrix then the 
eigenvector $v$ must have rational entries.  We begin with a few lemmas.

\begin{lem} \label{ratvalue}
If $A \in M_3(\Q)$ and $(a ~ b ~ 1)^T$ is an eigenvector of $A$ with rational 
eigenvalue $\lambda$ then $a$ and $b$ are rational.
\end{lem}

\begin{proof}
We can solve for $a$ and $b$ in terms of $\lambda$ and the entries of $A$.  Both $a$ 
and $b$ are rational functions of $\lambda$ and the $a_{ij}$ hence $a$ and $b$ are 
rational.
\end{proof}

\begin{lem} \label{power}
If $A \in GL_3(\R)$ has three linearly independent eigenvectors with corresponding 
eigenvalues $\lambda_1 > | \lambda_2 | \ge | \lambda_3 |$ then one of the columns of
$A^k$ converges to a multiple of the eigenvector $v_1$.
\end{lem}

\begin{proof}
Since $A$ is in the general linear group its column vectors must be linearly 
independent and therefore one of these columns, call it $u$, cannot lie in the span of
$v_2$ and $v_3$.  Hence we can write $u = c_1v_1 + c_2v_2 + c_3v_3$ with $c_1 \ne 0$.
\begin{eqnarray*}
A^k u &=& A^k (c_1v_1 + c_2v_2 + c_3v_3) \\
      &=& c_1 A^k v_1 + c_2 A^k v_2 + c_3 A^k v_3 \\
      &=& c_1 \lambda_1^k v_1 + c_2 \lambda_2^k v_2 + c_3 \lambda_3^k v_3 \\
      &=& \lambda_1^k( c_1 v_1 + c_2 (\lambda_2 / \lambda_1)^k v_2 + 
                       c_3 (\lambda_3 / \lambda_1)^k v_3 ) \\
\end{eqnarray*}
Because $\lambda_1 > | \lambda_2 | \ge | \lambda_3 |$ the limit of the column $A^ku$ 
approaches a multiple of $v_1$, namely $\lambda_1^k (c_1v_1)$.
\end{proof}

\begin{thm}
If $A \in GL_3(\Q)$ is stochastic and the columns of $A^k$ converge to a 
multiple of $(a ~ b ~ 1)^T$ then $a$ and $b$ are rational numbers.
\end{thm}

\begin{proof}
We know that $(a ~ b ~ 1)^T$ is an eigenvector of $A$.  Since $A$ is a stochastic 
matrix $A$ has an eigenvector with corresponding eigenvalue $1$, furthermore $1$ is 
the maximal eigenvalue \cite{Minc}.  Therefore the characteristic polynomial of $A$ is
$f(\lambda) = det(A - \lambda I) = 
(\lambda - 1)(\lambda - \lambda_2)(\lambda - \lambda_3)$.  To show that $a$ and $b$ 
are rational we consider a few cases.

Assume the $\lambda_i$ are distinct.  If the eigenvalues of $A$ satisfy
$1 > |\lambda_2| \ge |\lambda_3|$ then we can apply the previous lemma, 
concluding  $(a ~ b ~ 1)$ is the eigenvector with eigenvalue $1$.  Using
lemma \ref{ratvalue}, we infer that $a$ and $b$ are rational.  Now suppose 
that $| \lambda_2 | \ge 1$.  Then $|\lambda_2| = 1$ because $1$ is the maximal 
eigenvalue.  Because the $\lambda_i$ are distinct either $\lambda_2 = -1$ or 
$\lambda_2$ and $\lambda_3$ are complex conjugates.  In both cases, after 
diagonalizing $A$ we see that the columns of $A^k$ will not converge, contradicting 
our assumption.

If the $\lambda_i$ are not distinct then either $1 = \lambda_2$ or 
$1 = \lambda_2 = \lambda_3$.  Since $f(x)$ has rational coefficients $\lambda_3$ must
be rational.  Since all three eigenvalues are rational the associated eigenvectors
$(a ~ b ~ 1)^T$ must have rational entries.
\end{proof}

\begin{thm}\label{rational}
If $\{ \tilde{a}_k(i_k) \}$ is an eventually periodic weighted Bary sequence that 
converges to the point $(\alpha,\beta) \in \Tri$ then $\alpha$ and $\beta$ are 
rational.
\end{thm}

\begin{proof}
Again let $A$ denote the product of the initial
non-periodic matrices and $B$ be the product of the periodic part.  Then a subsequence
of the Bary partition triangles about the point $(\alpha,\beta)$ is
\[
A,AB,AB^2,AB^3,...
\]

Since $\{ \tilde{a}_k(i_k) \}$ converges to $(\alpha,\beta)$ we know that each column 
of the matrices $B,B^2,B^3,...$ must converge to a multiple of 
$A^{-1}(\alpha ~ \beta ~ 1)^T$.  Because $B$ is a product of stochastic transformation
matrices $B$ is also stochastic.  By the previous theorem, 
$A^{-1}(\alpha ~ \beta ~ 1)$ has rational entries.  Since everything in sight is 
rational $\alpha$ and $\beta$ must be rational.
\end{proof}

\section{The Farey-Bary Analog of Singularness}

The Minkowski question-mark function is a naturally occurring singular function.  It 
is strictly increasing, continuous, and has derivative zero almost everywhere.  Since 
$?(x)$ is a monotonic,
real-valued function acting on the unit interval we know that its derivative must
exist almost everywhere \cite{Royden}.  With this knowledge, 
proving $?'(x)=0$ a.e. was reduced to proving that if the following limit
existed and were finite then
\[
\mylim{n}{\infty} \frac{ \text{length of interval in the range} }
                       { \text{length of interval in the domain} } = 0,
\]
almost everywhere.  For the weighted Farey-Bary map we will show a stronger analog.  
Specifically, we will prove that
\[
\myliminf{n}{\infty} \frac{ \text{area of subtriangle in the range} }
                          { \text{area of subtriangle in the domain} } = 0,
\]
almost everywhere.

\subsection{Areas of the Farey and Bary Subtriangles}

We will frequently use the following theorem to calculate
the areas of the Farey and Bary subtriangles.

\begin{thm} \cite{fareybary} \label{area}
If $T$ is a triangle with vertices $(p_1/r_1,q_1/r_1),(p_2/r_2,q_2/r_2),$ and 
$(p_3/r_3,q_3/r_3)$ then $T$ can be represented by the matrix
\[
M = 
\left(
\begin{array}{ccc}
p_1 & p_2 & p_3 \\
q_1 & q_2 & q_3 \\
r_1 & r_2 & r_3 
\end{array}
\right),
\]
and
\[
area(T) = \frac{1}{2} \cdot \frac{|det(M)|}{r_1r_2r_3}.
\]
\end{thm}

\begin{figure}[hbtp]
\centering
\includegraphics[width=2.5in]{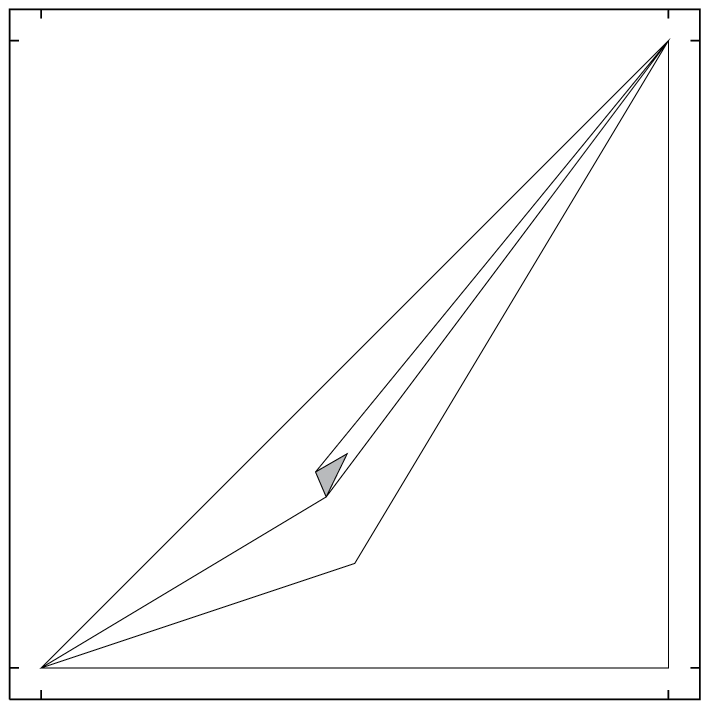}
\caption{$\Tri \{2(III),1(II),1(II)\}$.}
\label{wsequencefig}
\end{figure}

Before we can calculate the areas of the Farey and Bary subtriangles we need to 
introduce some new notation.  Given a finite Farey sequence 
$\{ a_1(i_1), ... , a_k(i_k) \}$ we define
\[
\Tri_k = \{ (x,y) : \{a_1(i_1),...a_k(i_k)\} 
\text{ are the first $k$ terms in the Farey sequence} \}.
\]
In figure \ref{wsequencefig} we see the set of points
$\Tri \{2(III),1(II),1(II)\}$.  We can think of this set by starting with the 
initial triangle $\Tri$ we perform a type $i_1$ transformation then $a_1 - 1$ type I 
operations continuing down the sequence we find the subtriangle that contains all 
the points whose Farey sequences begin with $\{ a_1(i_1),...,a_k(i_k \}$.
So, the matrix representation $M_{s_k}$ of $\Tri_k$ is a 
product of the initial matrix $M_0$ and $s_k = a_1 + ... + a_k$ transformation 
matrices.

\begin{cor}
Given a finite Farey sequence $\{ a_1(i_i), ... , a_k(i_k) \}$
\[
area( \Tri_k ) = \frac{1}{2} \cdot \frac{|det(M_{s_k})|}{r_1(k)r_2(k)r_3(k)}.
\]
\end{cor}

The determinant of $M_{s_k}$ will be a product of $s_k$ weights $m_1,m_2,$ and $m_3$,
dependent upon what transformation matrix was used.  Referring back to our Farey 
transformation matrices, taking a type I move multiplies the determinant of $M_{s_k}$
by $m_3$, type II contributes a factor of $m_1$, and type III provides an $m_2$ term.
The values for $r_i(k)$ are more difficult to calculate, we will deal with them in our
proof of singularness

\begin{cor}
Given a finite Bary sequence $\{ \tilde{a}_1(i_1), ... , \tilde{a}_k(i_k) \}$
\[
area( \tilde{\Tri}_k ) = \frac{1}{2} \cdot 
                         \frac{|det(\tilde{M}_{s_k})|}{1\cdot1\cdot1}.
\]
\end{cor}
           
Notice that the bottom row of the $\tilde{M}_{s_k}$ is filled with ones.  This is a
consequence of the definition of the weighted Bary sum.  The determinant of 
$\tilde{M}_{s_k}$ is a product of $w_1,w_2,$ and $w_3$, where 
$w_i = m_i/(m_1+m_2+m_3)$.  Looking back at the Bary transformation matrices, we see 
that taking a type I move multiplies the determinant by $w_3$, type II gives a factor 
of $w_1$, and type III supplies a term of $w_2$.

We are interested in calculating the ratio of the area of a Bary subtriangle (a 
subtriangle in the range of $\delta$) and the corresponding Farey subtriangle (in the 
domain)
\[
\frac{area( \tilde{\Tri} \{a_1(i_1),...,a_k(i_k)\} )}
     {area( \Tri\{a_1(i_1),...,a_k(i_k)\} )} 
=
\frac{ |det(\tilde{M}_{s_k})| }{ |det(M_{s_k})| } \cdot r_1(k)r_2(k)r_3(k).
\]
To simplify the ratio of determinants consider what happens when we make a type I 
move.  The determinant of $\tilde{M}_{s_k}$ is multiplied by $w_3 = m_3/(m_1+m_2+m_3)$
and the determinant of $M_{s_k}$ is multiplied by $m_3$; the net effect is to multiply
the ratio by $1/(m_1+m_2+m_3)$.  The same occurs when taking type II or III 
moves, therefore
\[
\frac{area( \tilde{\Tri}\{a_1(i_1),...,a_k(i_k)\} )}
     {area( \Tri\{a_1(i_1),...,a_k(i_k)\} )} 
=
\frac{ r_1(k)r_2(k)r_3(k) }{ (m_1+m_2+m_3)^{s_k} }.
\]
Much of this chapter will be devoted to showing when this ratio goes to zero.

\subsection{The Failure of the Farey Partition when $m_3 > 1$}

In this section, we show that if $m_3>1$ then the set of points that have infinite
Farey sequences has measure zero.  This is an interesting section in that it reveals
an asymmetry of the Farey-Bary map.  Now we set the stage for a very
useful lemma.

Given a finite Farey sequence $\{ a_1(i_1),...,a_k(i_k) \}$, let 
$T = \Tri \{ a_1(i_1),...,a_k(i_k) \}$ and
\[
M = 
(v_1 ~ v_2 ~ v_3) = 
\left(
\begin{array}{ccc}
p_1 & p_2 & p_3 \\
q_1 & q_2 & q_3 \\
r_1 & r_2 & r_3 \\
\end{array}
\right)
\]
be the matrix representation of $T$.  For any integer $L \ge 1$, define $T_L(i)$ to 
be the result of performing a type $i$ move and then $L-1$ type I moves on $T$.  Thus

\begin{eqnarray*}
T_L(I) &=& \langle v_1, v_2, (1+...+m_3^{L-2})(m_1v_1 \fsum m_2v_2) \fsum
                             m_3^{L-1}(m_1v_1 \fsum m_2v_2 \fsum m_3v_3) \rangle, \\
T_L(II) &=& \langle v_2, v_3, (1+...+m_3^{L-2})(m_1v_2 \fsum m_2v_3) \fsum 
                             m_3^{L-1}(m_1v_1 \fsum m_2v_2 \fsum m_3v_3) \rangle, \\
T_L(III) &=& \langle v_1, v_3, (1+...+m_3^{L-2})(m_1v_1 \fsum m_2v_3) \fsum 
                             m_3^{L-1}(m_1v_1 \fsum m_2v_2 \fsum m_3v_3) \rangle.
\end{eqnarray*}

Define $T_L = T - T_L(I) - T_L(II) - T_L(III)$.  $T_L$ contains all the points 
with $\{ a_1(i_1),...,a_k(i_k) \}$ as the first $k$ terms of their Farey sequence and 
$a_{k+1} < L$.

\begin{figure}[hbtp]
\centering
\includegraphics[width=\textwidth]{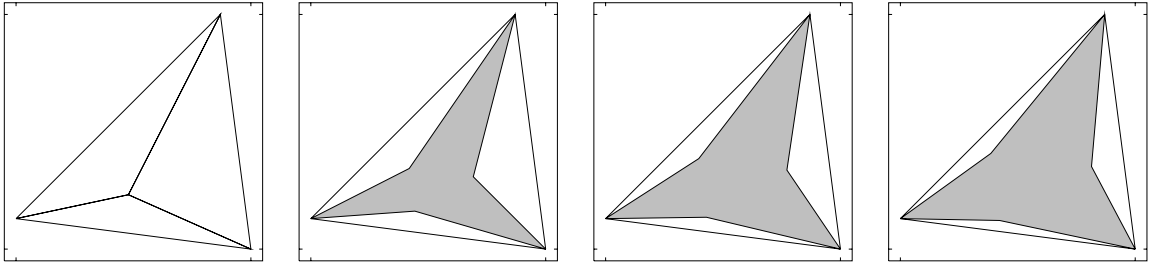}
\caption{$T_1,T_2,T_3$ and $T_4$.}
\end{figure}

\begin{lem}\label{T_L}
If $m_3 = 1$ and $L \ge 1$ then 
\[
\frac{area(T_L)}{area(T)} \le \frac{L-1}{L-1+\frac{1}{m_1+m_2}}.
\]
However, if $m_3 > 1$ and $L \ge 1$ then
\[
\frac{area(T_L)}{area(T)} \le \frac{1}{1+\frac{1}{m_1+m_2}}.
\]
\end{lem}

\begin{proof}
For ease of notation, we set $x = r_1, y = r_2, z = r_3.$  We know from theorem 
\ref{area} that
\[
2 \cdot area(T) = \frac{|det(M)|}{xyz}.
\]
To bound the area of $T_L = area(T) - area(T_L(I)) - area(T_L(II)) - area(T_L(III))$, 
we begin by calculating the areas of the various $T_L(i)$ from their matrix 
representations
\begin{eqnarray*}
2 \cdot area(T_L) = \frac{|det(M)|}{xyz} &-& \frac{|det(M)|m_3m_3^{L-1}}{[(1 + ... + m_3^{L-2})(m_1x + m_2y) + m_3^{L-1}(m_1x + m_2y + m_3z)]xy} \\
                                         &-& \frac{|det(M)|m_1m_3^{L-1}}{[(1 + ... + m_3^{L-2})(m_1y + m_2z) + m_3^{L-1}(m_1x + m_2y + m_3z)]yz} \\
                                         &-& \frac{|det(M)|m_2m_3^{L-1}}{[(1 + ... + m_3^{L-2})(m_1x + m_2z) + m_3^{L-1}(m_1x + m_2y + m_3z)]xz}.
\end{eqnarray*}
Now we pull out a common factor of $2 \cdot area(T)$
\begin{eqnarray*}
2 \cdot area(T_L) = \frac{|det(M)|}{xyz} \Bigg[ 1 &-& \frac{m_3^{L-1}m_3z}{(1 + ... + m_3^{L-2})(m_1x + m_2y) + m_3^{L-1}(m_1x + m_2y + m_3z)} \\
                                             &-& \frac{m_3^{L-1}m_1x}{(1 + ... + m_3^{L-2})(m_1y + m_2z) + m_3^{L-1}(m_1x + m_2y + m_3z)} \\
                                             &-& \frac{m_3^{L-1}m_2y}{(1 + ... + m_3^{L-2})(m_1x + m_2z) + m_3^{L-1}(m_1x + m_2y + m_3z)} \Bigg].
\end{eqnarray*}
We make the right hand side larger by using the largest denominator for all three 
fractions
\begin{eqnarray*}
2 \cdot area(T_L) &\le& \frac{|det(M)|}{xyz} \left[ 1 - \frac{m_3^{L-1}(m_1x + m_2y + m_3z)}{[(1 + ... + m_3^{L-2})(m_1y + m_2z) + m_3^{L-1}(m_1x + m_2y + m_3z)]} \right] \\
                  &=&   \frac{|det(M)|}{xyz} \frac{(1 + ... + m_3^{L-2})(m_1y + m_2z)}{(1 + ... + m_3^{L-2})(m_1y + m_2z) + m_3^{L-1}(m_1x + m_2y + m_3z)}.
\end{eqnarray*}

We now have two cases to consider.  First, suppose that $m_3 = 1$ then
\begin{eqnarray*}
\frac{area(T_L)}{area(T)} &\le& \frac{(L-1)(m_1y + m_2z)}{(L-1)(m_1y + m_2z) + (m_1x + m_2y + m_3z)} \\
                          &=&   \frac{L-1}{L-1 + \frac{m_1x + m_2y + m_3z}{m_1y + m_2z}} \\
                          &\le& \frac{L-1}{L-1 + \frac{m_3z}{m_1z + m_2z}} \\
                          &=&   \frac{L-1}{L-1 + \frac{1}{m_1 + m_2}}.
\end{eqnarray*}
To prove our second bound, suppose that $m_3 > 1$ then we know for all $L \ge 2$ that
\[
1 + m_3 + ... + m_3^{L-2} = \frac{m_3^{L-1} - 1}{m_3 - 1}.
\]
Substituting this back into our intermediate result we have
\begin{eqnarray*}
\frac{area(T_L)}{area(T)} &\le& \frac{\frac{m_3^{L-1}-1}{m_3-1}(m_1y + m_2z)}{\frac{m_3^{L-1}-1}{m_3-1}(m_1y + m_2z) + m_3^{L-1}(m_1x + m_2y + m_3z)} \\
                          &=&   \frac{(m_3^{L-1}-1)(m_1y + m_2z)}{(m_3^{L-1}-1)(m_1y + m_2z) + m_3^{L-1}(m_3 - 1)(m_1x + m_2y + m_3z)} \\
                          &\le& \frac{(m_3^{L-1}-1)(m_1y + m_2z)}{(m_3^{L-1}-1)(m_1y + m_2z) + m_3^{L-1}(m_1x + m_2y + m_3z)} \\
                          &\le& \frac{(m_3^{L-1}-1)(m_1y + m_2z)}{(m_3^{L-1}-1)(m_1y + m_2z) + (m_3^{L-1}-1)(m_1x + m_2y + m_3z)} \\
                          &=&   \frac{(m_1y+m_2z)}{(m_1y+m_2z)+(m_1x+m_2y+m_3z)} \\
                          &=&   \frac{1}{1+\frac{m_1x+m_2y+m_3z}{m_1y+m_2z}} \\
                          &\le& \frac{1}{1+\frac{m_3z}{(m_1+m_2)z}} \\
                          &=&   \frac{1}{1+\frac{m_3}{m_1+m_2}} \\
                          &\le& \frac{1}{1+\frac{1}{m_1+m_2}}.
\end{eqnarray*}
\end{proof}

\begin{figure}[hbtp]
\centering
\includegraphics[height=2.5in]{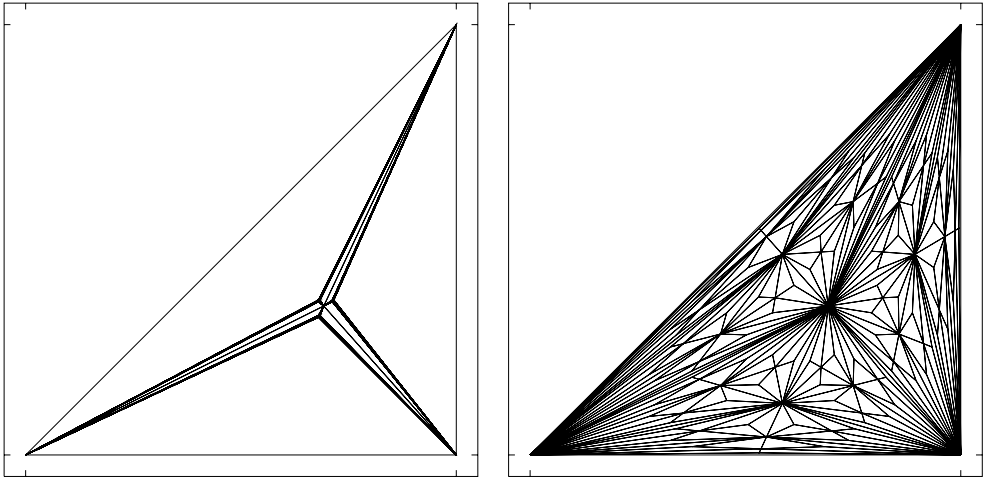}
\caption{The failure of the Farey partition when $m_3 > 1$.}
\end{figure}

Since the area of $T_L$ never reaches the full area of $T$ it appears that the 
Farey partition does not fully partition the triangle when $m_3 > 1$.  We now state
and prove what we mean by the failure of the Farey partition.

\begin{thm}
If $m_3 > 1$ then the set of points $(\alpha,\beta) \in \Tri$ that have infinite Farey
sequences has measure zero.
\end{thm}

\begin{proof}
Let $M$ be the set of points that have infinite Farey sequences.  For 
$(\alpha,\beta) \in M$ to have an infinite Farey sequence $\{a_1(i_1),a_2(i_2),...\}$
each $a_k$ must be finite.

Set
\[
M(1) = \{ (\alpha,\beta) \in \Tri : a_1 < \infty \},
\]
and in general
\[
M(k) = \{ (\alpha,\beta) \in M(k-1) : a_k < \infty \}.
\]
We have just defined a nested sequence of sets with
\[
M = \bigcap_{k=1}^\infty M(k).
\]

From the previous lemma we have
\[
measure(M(k)) \le \frac{1}{1+\frac{1}{m_1+m_2}}measure(M(k-1)).
\]
Iterating this inequality gives a bound on the measure of $M$
\[
measure(M) \le \prod_{k=1}^\infty\frac{1}{1+\frac{1}{m_1+m_2}} = 0.
\]
\end{proof}

For the rest of the paper we will assume that $m_3 = 1$.

\subsection{The Definition of the Weighted Farey-Bary Map}\label{definition}

Our first task in this section is to show that the Farey partition does not fail when 
$m_3 = 1$.  The second is to define the weighted Farey-Bary map.

\begin{thm}
If $m_3 = 1$ and $v \in \Tri_F$ is an interior point then the Farey sequence 
associated to $v$ is infinite.
\end{thm}

\begin{proof}
Assume the Farey sequence for $v$ is finite, then the last term in the sequence must be 
infinite.  Let $\langle v_1,v_2,v_3 \rangle$ be the triangle that contains $v$ before
we perform the infinite sequence of type I moves.  Then the final Farey subtriangles
containing $v$ are given by
\[
v_1, v_2, (L-1)(m_1v_1 \fsum m_2v_2) \fsum (m_1v_1 \fsum m_2v_2 \fsum m_3v_3),
\]
where $L$ is the number of type I operations performed.
This sequence converges to the line segment
\[
v_1, v_2, m_1v_1 \fsum m_2v_2.
\]
Thus $v$ must lie on this edge contradicting the fact that $v$ is an interior point.
\end{proof}

\begin{cor}
If $m_3=1$ then almost all points $(\alpha,\beta) \in \Tri$ have infinite Farey
sequences.
\end{cor}

\begin{thm}
If $v \in \Tri_F$ is an interior point then the sequence $\{ \delta_n(v) \}$ converges
to a point in $\Tri_B$.
\end{thm}

\begin{proof}
Let $v \in \Tri_F$ be an interior point.  For each $n$, $v$ must land on the interior
of one of the $3^n$ subtriangles of $P_n$.  Label this triangle $\Tri_n(v)$, and 
denote the corresponding subtriangle in $\tilde{P}_n$ by $\tilde{\Tri}_n(v)$.  We will
show that the vertices of $\tilde{\Tri}_n(v)$ converge to a point and therefore the
sequence $\{ \delta_n(v) \}$ must also converge.

A calculation like those we will use in lemma 3.2.2 shows that the distance from the
barycenter of a triangle to any of the vertices is at most $(m_1+m_2)/(m_1+m_2+m_3)$ 
times the length of the longest side of the triangle.  Because each interior point has
an infinite Farey sequence we must perform an infinite number of type II or III moves.
Consequently, at most one of the vertices in $\tilde{\Tri}_n(v)$ will remained fixed 
for infinitely many of the $n$.  This forces the lengths of the sides of 
$\tilde{\Tri}_n(v)$ to approach zero, hence the sequence $\delta_n(v)$ must also 
converge.
\end{proof}

\begin{defn}
Define the weighted Farey-Bary map $\delta : \Tri_F \To \Tri_B$ by setting
\[
\delta(v) = \mylim{n}{\infty} \delta_n(v)
\]
when $v$ is an interior point and
\[
\delta(v) = \delta_n(v)
\]
when $v$ lies on an edge of one of the partition triangles in $P_n$ for some 
$n$.
\end{defn}

While the weighted Farey-Bary map is not continuous at all points, it is continuous 
almost everywhere.

\begin{thm}
$\delta$ is continuous at all interior points.
\end{thm}

\begin{proof}
Let $v$ be an interior point, and $\{ v_k \}$ be a sequence of points in $\Tri_F$ that
converges to $v$.  Given any small ball about $\delta(v)$ there exists an $n$ such 
that $\tilde{\Tri}_n(v)$ is contained within this ball because the lengths of the 
sides of $\tilde{\Tri}_n(v)$ approach zero.  The inverse image of $\tilde{\Tri}_n(v)$ 
is $\Tri_n(v)$, which contains $v$.  Since $\{ v_k \}$ converges to $v$ the $v_k$ will
eventually land within $\Tri_n(v)$.  Hence the $\delta(v_k)$ will eventually be in 
$\tilde{\Tri}_n(v)$, inside our small ball about $\delta(v)$.  Thus 
$\{ \delta(v_k) \}$ converges to $\delta(v)$, concluding $\delta$ is continuous at 
$v$.
\end{proof}

\subsection{Almost everywhere $\limsup \frac{a_1+...+a_n}{n} = \infty$}

Our immediate goal is to show that the sum of the terms in a Farey sequence
grows faster than the length of the sequence.  This fact is crucial to our proof of 
singularness.  

\begin{thm}
The set of $(\alpha,\beta) \in \Tri$ for which
\[
\mylimsup{n}{\infty} \frac{a_1 + ... + a_n}{n} < \infty
\]
has measure zero.
\end{thm}

\begin{proof}
For each positive integer $N$, set
\[
M_N = \{(\alpha,\beta) \in \Tri : \forall n \ge 1, \frac{a_1 + ... + a_n}{n} \le N \}.
\]
Since the union of all the $M_N$ is the set we want to show has measure zero if we can
show that the $measure(M_N) = 0$ then we will be done.

Now, $\frac{a_1 + ... + a_n}{n} \le N$ if and only if $a_1 + ... + a_n \le nN$.  
Certainly, we can say for all $n$ that
\[
a_n \le nN + 1.
\]
Thus 
\[
M_N \subset \tilde{M}_N = \{ (\alpha,\beta) \in \Tri:\forall n \ge 1, a_n \le nN+1 \}.
\]
We have reduced this proof to showing that the $measure(\tilde{M}_N) = 0$.

Set
\[
\tilde{M}_N(1) = \{ (\alpha,\beta) \in \Tri : a_1 \le N + 1\},
\]
and in general
\[
\tilde{M}_N(k) = \{ (\alpha,\beta) \in \tilde{M}_N(k-1) : a_k \le kN + 1\}.
\]
Then we have a nested sequence of sets with
\[
\tilde{M}_N = \bigcap_{k=1}^\infty \tilde{M}_N(k).
\]

But this puts us into the language of lemma \ref{T_L}.  Letting $L = kN + 1$, 
we can conclude that
\[
measure(\tilde{M}_N(k)) \le \frac{kN}{kN + c}measure(\tilde{M}_N(k-1))
\]
where $c = \frac{1}{m_1+m_2}$, and hence
\[
measure(\tilde{M}_N(k)) \le \prod_{k=1}^\infty \frac{kN}{kN + c}.
\]
We must show that this infinite product is zero, which is equivalent to showing that 
its reciprocal
\[
\prod_{k=1}^\infty \frac{kN + c}{kN} = \prod_{k=1}^\infty (1 + \frac{c}{kN}) = \infty.
\] 
Taking logarithms, this is the same as showing that the series
\[
\sum_{k=1}^\infty log(1 + \frac{c}{kN}) = \infty.
\]
We know this diverges since, for large enough $k$, we have
\begin{eqnarray*}
log(1+\frac{c}{kN}) &=& log(1+\frac{1}{kN(m_1+m_2)}) \\
                    &\ge& \frac{1}{2kN(m_1+m_2)},
\end{eqnarray*}
which is a multiple of the harmonic series.
\end{proof}

\subsection{Almost everywhere $\liminf (area(\tilde{\Tri}_n)/area(\Tri_n)) = 0$}

We have reached the climax of this chapter.  We are about to prove a direct 
generalization that the Minkowski question-mark function is singular, specifically
\[
\myliminf{n}{\infty}\frac{\text{area of a subtriangle in the range}}
                         {\text{area of a subtriangle in the domain}} = 0,
\]
almost everywhere.  Intuitively, we are showing that at almost all points
\[
\myliminf{n}{\infty} det( \text{Jacobian of } \delta_n ) = 0.
\]
The proof of the following theorem was inspired by Viader, Paradis, and Bibiloni's
work in \cite{NewLight}.

\begin{thm}
For any point $(\alpha, \beta) \in \Tri$, off a set of measure zero,
\[
\myliminf{n}{\infty} \frac{area(\tilde{\Tri}\{a_1(i_1), ... , a_n(i_n)\})}
                          {area(\Tri\{a_1(i_1),...,a_n(i_n)\})} = 0.
\]
\end{thm}

\begin{proof}
Letting $s_n = a_1 + ... + a_n$, we have already shown that
\[
\frac{area( \Tri\{\tilde{a}_1(i_1),...,\tilde{a}_n(i_n)\} )}
     {area( \Tri\{a_1(i_1),...,a_n(i_n)\} )} 
=
\frac{ r_1(n)r_2(n)r_3(n) }{ (m_1+m_2+m_3)^{s_n} }.
\]
We need to bound the size of the $r_i(n)$.  Returning to our recurrence relations, we
have three possibilities for the value of $r_3(n)$ 
\begin{enumerate}
\centering
\item $(a_n-1)(m_1r_1(n-1) + m_2r_2(n-1)) + (m_1r_1(n-1) + m_2r_2(n-1) + m_3r_3(n-1))$
\item $(a_n-1)(m_1r_2(n-1) + m_2r_3(n-1)) + (m_1r_1(n-1) + m_2r_2(n-1) + m_3r_3(n-1))$
\item $(a_n-1)(m_1r_1(n-1) + m_2r_3(n-1)) + (m_1r_1(n-1) + m_2r_2(n-1) + m_3r_3(n-1))$
\end{enumerate}
Therefore
\[
r_3(n) \le a_n(m_1+m_2+m_3)r_3(n-1),
\]
iterating this inequality gives us
\[
r_3(n) \le (m_1+m_2+m_3)^n \cdot \prod_{i=1}^n a_i.
\]

By convention we always have $r_1(n) \le r_2(n) \le r_3(n)$, thus
\[
\frac{area(\tilde{\Tri}\{a_1(i_1),...,a_n(i_n)\})}
     {area(\Tri\{a_1(i_1),...,a_n(i_n)\})} 
\le 
\frac{(m_1+m_2+m_3)^{3n}\prod a_i^3}{(m_1+m_2+m_3)^{s_n}}.
\]
By the arithmetic-geometric mean we know
\[
\prod_{i=1}^n a_i \le \left( \frac{a_1 + ... + a_n}{n} \right)^n 
= \left( \frac{s_n}{n} \right)^n.
\]
Applying this to the previous inequality, we have
\begin{eqnarray*}
\frac{area(\tilde{\Tri}\{a_1(i_1),...,a_n(i_n)\})}
     {area(\Tri\{a_1(i_1),...,a_n(i_n)\})} 
&\le& 
\frac{(m_1+m_2+m_3)^{3n}(\frac{s_n}{n})^{3n}}{(m_1+m_2+m_3)^{s_n}} \\
&=&
\left( \frac{(m_1+m_2+m_3)^3(\frac{s_n}{n})^3}{(m_1+m_2+m_3)^{s_n/n}} \right)^n.
\end{eqnarray*}

From the previous theorem, we know that $s_n/n \To \infty$, almost everywhere.  Since 
the above denominator has a $(m_1+m_2+m_3)^{s_n/n}$ term while the numerator only has 
a $(s_n/n)^3$ term, the entire ratio must approach zero.
\end{proof}


\chapter{The Continuous Farey-Bary Map}

The major shortcoming of the weighted Farey-Bary map, as with the original Farey-Bary
map, is its discontinuities.  In 
this chapter, we develop a new partitioning scheme that will create a continuous
Farey-Bary map.  Some of our previous results follow through to this chapter.  For 
instance, the continuous Farey-Bary map will send a natural class of cubic irrational 
points to a natural class of rational points.  But, our proof of singularness does 
not carry over.  We conclude the chapter by proving a weaker analog of singularness.

\section{The Revised Farey and Bary Partitions}

Again we define two partitions of the triangle
\[
\Tri = \{ (x,y) : 1 \ge x \ge y \ge 0 \}.
\]
The first partition will yield the domain of our desired function, and the second 
will yield the range.

We define the revised Farey partition of the triangle by a sequence of partitions
$\{ P_n \}$ such that each $P_n$ will consist of $6^n$ subtriangles of $\Tri$ and each
$P_n$ will be a refinement of the previous $P_{n-1}$.  Let $P_0$ be the initial 
triangle $\Tri$.  
\begin{eqnarray*}
v_1 &=& (0,0,1), \\
v_2 &=& (1,0,1), \\
v_3 &=& (1,1,1). 
\end{eqnarray*}
Now we take all possible Farey sums of these three vertices producing four new points
\begin{eqnarray*}
&v_1 \fsum v_2 = (1,0,2) = (\frac{1}{2},\frac{0}{2}),& \\
&v_1 \fsum v_3 = (1,1,2) = (\frac{1}{2},\frac{1}{2}),& \\
&v_2 \fsum v_3 = (2,1,2) = (\frac{2}{2},\frac{1}{2}),& \\
&v_1 \fsum v_2 \fsum v_3 = (2,1,3) = (\frac{2}{3},\frac{1}{3}).&
\end{eqnarray*}
These four points refine $P_0$ by cutting the triangle into six new subtriangles.  
Thus determining the partition $P_1$.  

We now proceed inductively.  Suppose we have the partition $P_n$ that determines $6^n$
subtriangles of $\Tri$.  We split each of these triangles into six smaller 
subtriangles as follows.  Given a subtriangle of $P_n$ with vertices 
$v_1,v_2,v_3$, we compute all four possible Farey sums of the three vertices producing
four new points that break the triangle into six new subtriangles.  Calculating the 
Farey sums for each
subtriangle of $P_n$, gives us the next partition $P_{n+1}$ of $\Tri$.  We denote this
full partitioning of $\Tri$ by $\Tri_F$ and call it the Farey partitioning giving us 
the domain of the continuous Farey-Bary map.

For the revised Bary partition we define a similar sequence of partitions 
$\{ \tilde{P}_n \}$ of $\Tri$ such that each $\tilde{P}_n$ will consist of $6^n$ 
triangles of $\Tri$ and each $\tilde{P}_n$ will be a refinement of the previous 
$\tilde{P}_{n-1}$.  The only difference here, is that we will take all possible Bary
sums, and so we call this partitioning of $\Tri$ the Bary partitioning and denote it 
$\Tri_B$.

\begin{figure}[p]
\centering
\includegraphics[height=.9\textheight]{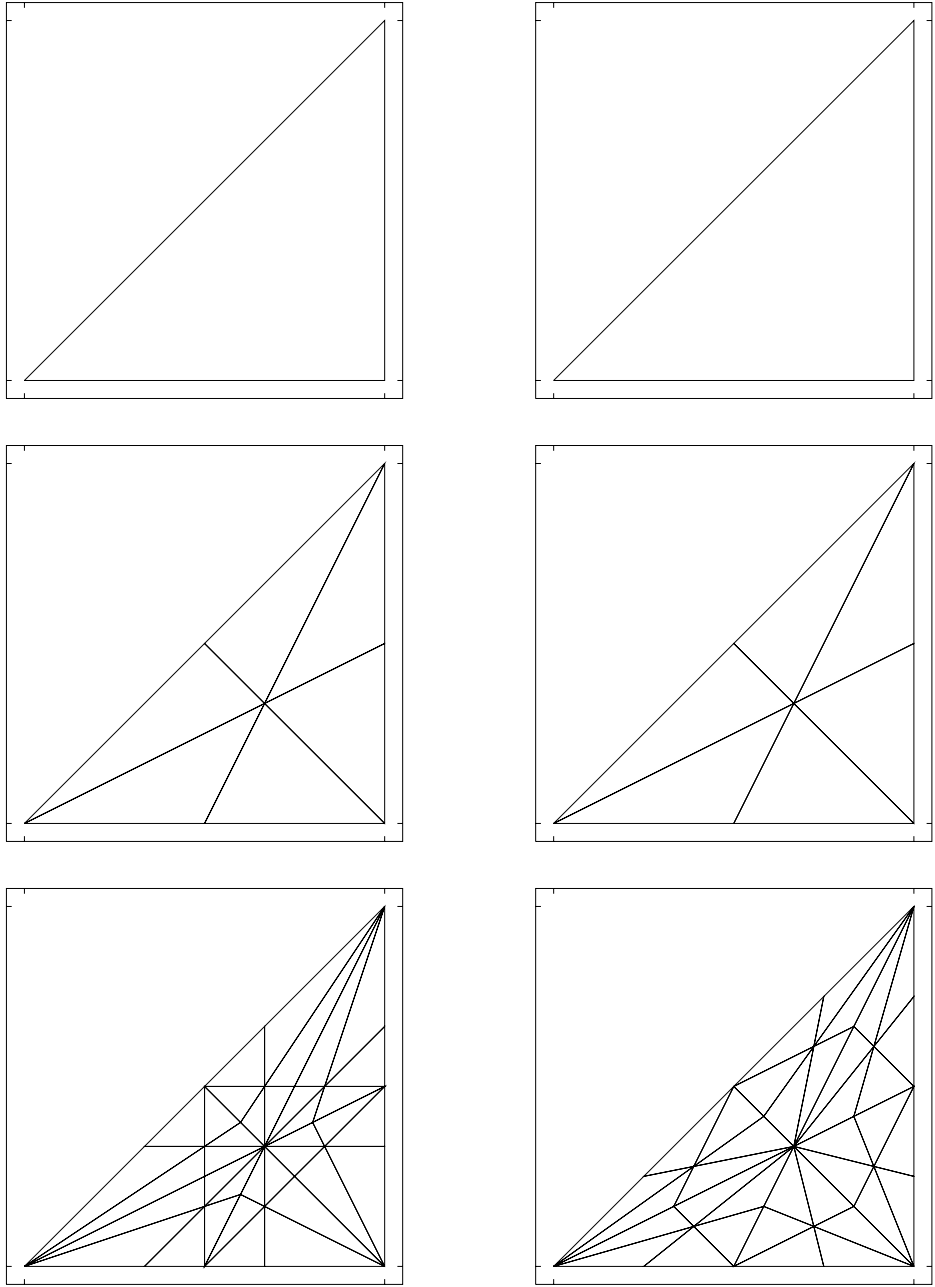}
\caption[The progression of the revised Farey and Bary partitions.]
        {The first three iterations of the revised Farey and Bary partitions.}
\label{revill}
\end{figure}

\begin{figure}[p]
\centering
\includegraphics[height=.9\textheight]{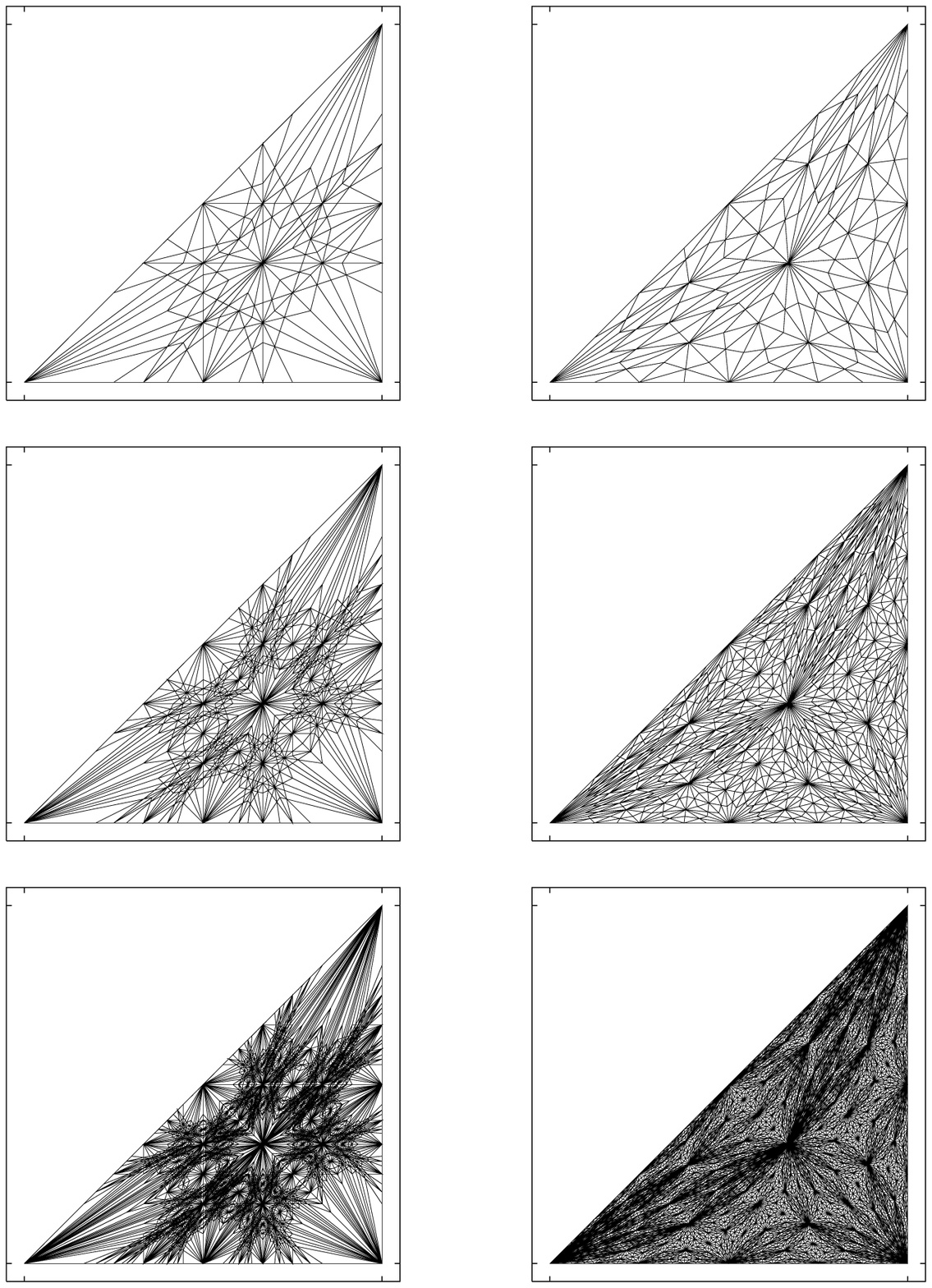}
\caption{The revised Farey and Bary partitions continued.}
\label{revill2}
\end{figure}

We illustrate the first six iterations of these partitions in figures \ref{revill}
and \ref{revill2}.  We now define the continuous Farey-Bary map by mapping the Farey
partition onto the Bary partition.

\section{The Continuous Farey-Bary Map}

Recall that we denote an arbitrary subtriangle in the $n^{th}$ stage of the Farey 
partitioning by $\langle v_1(n), v_2(n), v_3(n) \rangle$.  And for an arbitrary 
subtriangle in the Bary partition we use 
$\langle \tilde{v}_1(n), \tilde{v}_2(n), \tilde{v}_3(n) \rangle$.

\begin{defn}
Given a natural number $n$, define $\delta_n$ to send any vertex in the $n^{th}$ Farey
partition $P_n$ to the corresponding vertex in the $n^{th}$ Bary partition
$\tilde{P}_n$.  That is, define $\delta_n$ on any subtriangle 
$\langle v_1(n), v_2(n), v_3(n) \rangle$ of $P_n$ by
\[
\delta(v_i(n)) = \tilde{v}_i(n).
\]
Finally, for any point $(x,y)$ in the subtriangle with vertices 
$\langle v_1(n), v_2(n), v_3(n) \rangle$, set
\[
\delta(x,y) = \alpha \tilde{v}_1(n) + \beta \tilde{v}_2(n) + \gamma \tilde{v}_3(n),
\]
where
\[
(x,y) = \alpha v_1(n) + \beta v_2(n) + \gamma v_3(n).
\]
\end{defn}

\subsection{Continuity}

We will now prove that the continuous Farey-Bary map is continuous and thus deserves
its name.  This is the first and foremost improvement on the original Farey-Bary map.

\begin{lem}\label{sidelength}
If $\langle \tilde{v}_1(n), \tilde{v}_2(n), \tilde{v}_3(n) \rangle$ is a subtriangle 
of the $n^{th}$ Bary partition $\tilde{P}_n$ then the length of its longest side is 
less than or equal to $(2/3)^{n-1}$.
\end{lem}

\begin{proof}
For the initial Bary partition $\tilde{P}_0 = \Tri$, the longest side of this triangle
has length $\sqrt{2} \le (2/3)^{-1}$.  Now consider what happens when we go from the 
$(n-1)^{st}$ partition to the $n^{th}$.  Let 
$\tilde{\Tri}_{n-1} = \langle \tilde{v}_1,\tilde{v}_2,\tilde{v}_3 \rangle$ 
be an arbitrary subtriangle of the $(n-1)^{st}$ Bary partition.  By assumption the 
longest
side of $\tilde{\Tri}_{n-1}$ has length $l \le (2/3)^{n-2}$  After subdividing 
$\tilde{\Tri}_{n-1}$ into six subtriangles, we have twelve new sides to bound.  For 
the six exterior sides we have
\begin{eqnarray*}
d(v_i, v_i \bsum v_j) &=& \sqrt{ (x_i/2 - x_j/2)^2 + (y_i/2 - y_j/2)^2 } \\
                      &=& \frac{1}{2} \sqrt{ (x_i - x_j)^2 + (y_i - y_j)^2 } \\
                      &=& \frac{1}{2} d(v_i,v_j) \\
                      &\le& (1/2) \cdot l \\
                      &\le& (2/3)^{n-1}.
\end{eqnarray*}
For the three longer interior edges we have
\begin{eqnarray*}
d(v_i, v_i \bsum v_j \bsum v_k) 
        &=& \sqrt{ (2x_i/3 - x_j/3 - x_k/3)^2 + (2y_i/3 - y_j/3 - y_k/3)^2 } \\
        &=& \frac{2}{3} \sqrt{ (x_i - x_j/2 - x_k/2)^2 + (y_i - y_j/2 - y_k/2)^2 } \\
        &=& \frac{2}{3} d(v_i, v_j \bsum v_k) \\
        &\le& (2/3) \cdot l \\
        &\le& (2/3)^{n-1}.
\end{eqnarray*}
And for the three shorter interior edges we have
\begin{eqnarray*}
d(v_i \bsum v_j, v_i \bsum v_j \bsum v_k) &=& \frac{1}{3} d(v_k, v_i \bsum v_j) \\
                                          &\le& (2/3)^{n-1}.        
\end{eqnarray*}
Our bound holds in all cases.
\end{proof}

\begin{thm}
The sequence of functions $\{ \delta_n \}$ is uniformly convergent.
\end{thm}

\begin{proof}
We want to show that for every $\epsilon > 0$ there exists an $N$ such that for all 
$v \in \Tri$ and $m,n \ge N$ we have $d( \delta_m(v), \delta_n(v) ) \le \epsilon$.
From the previous lemma we know how to find an $N$ such that the length of the sides
of the subtriangles in $\tilde{P}_N$ are less than or equal to $\epsilon$.  Since 
$\delta_m(v)$ and $\delta_n(v)$ must land in the same subtriangle of $\tilde{P}_n$ we
are guaranteed that $d( \delta_m(v), \delta_n(v) ) \le \epsilon$
\end{proof}

Since the sequence of functions $\{ \delta_n \}$ is uniformly convergent we know that
its limit exists.

\begin{defn}
Define the continuous Farey-Bary map, $\delta : \Tri_F \To \Tri_B$, to be the limit of
the sequence $\{ \delta_n \}$.
\end{defn}

\begin{cor}
The continuous Farey-Bary map is continuous.
\end{cor}

\begin{proof}
Since each $\delta_n$ is continuous and because the sequence converges uniformly we 
know that the limit function must also be continuous.
\end{proof}

\subsection{Restricted to a Single Edge $\delta$ Acts Like $?(x)$}

One reason why we believe the continuous Farey-Bary map to be the natural extension 
of the Minkowski question-mark function is that when we look at $\delta$ on an edge of
one of the partition subtriangles it acts like $?(x)$.  In this section we
will add a little formality to this notion.

But first, consider what $\delta$ does on the X-axis.  In the first few stages of the 
Farey partition the base of $\Tri$ is broken up as follows.
\begin{eqnarray*}
& I_0 = \left[ (\frac{0}{1}, 0), (\frac{1}{1}, 0) \right] & \\
& I_1 = \left[ (\frac{0}{1}, 0), (\frac{1}{2}, 0), (\frac{1}{1}, 0) \right] & \\
& I_2 = \left[ (\frac{0}{1}, 0), (\frac{1}{3}, 0), (\frac{1}{2}, 0), (\frac{2}{3}, 0),
               (\frac{1}{1}, 0) \right] & \\
\end{eqnarray*}
In the Bary partition, we have
\begin{eqnarray*}
& I_0 = \left[ (\frac{0}{1}, 0), (\frac{1}{1}, 0) \right] & \\
& I_1 = \left[ (\frac{0}{1}, 0), (\frac{1}{2}, 0), (\frac{1}{1}, 0) \right] & \\
& I_2 = \left[ (\frac{0}{1}, 0), (\frac{1}{4}, 0), (\frac{1}{2}, 0), (\frac{3}{4}, 0),
               (\frac{1}{1}, 0) \right] & \\
\end{eqnarray*}
We can see that $\delta(x,0) = (?(x), 0)$.  Although we will not have as clean a
formula for all edges we will see a clear relation between the two functions.  First,
we prove a well known result.

\begin{lem}
If $gcd(a,b) = 1$ then there exist consecutive Farey fractions with $a$ and $b$ as 
their denominators.
\end{lem}

\begin{proof}
We will prove this lemma inductively.  Our base case is when $a=b=1$, which is 
satisfied by the Farey fractions $0/1$ and $1/1$.  Without loss of generality suppose 
that $a<b$, and assume that this lemma holds for all pairs of integers less than $a$.
Using the division algorithm we know
\[
b = qa + r,
\]
where $q > 0$ and $0 < r < a$.  So, $r$ and $a-r$ must be consecutive Farey fractions
by our inductive hypothesis.  At some stage we have
\begin{eqnarray*}
I_k: & \left[ ..., \frac{n_1}{r}, \frac{n_2}{a-r}, ... \right] & \\
I_{k+1}: & \left[ ..., \frac{n_1}{r}, \frac{n_1+n_2}{a}, \frac{n_2}{a-r}, ... \right] & \\
I_{k+2}: & \left[ ..., \frac{n_1}{r}, \frac{2n_1+n_2}{a+r}, \frac{n_1+n_2}{a}, ... \right] & \\
I_{k+3}: & \left[ ..., \frac{2n_1+n_2}{a+r}, \frac{3n_1+2n_2}{2a+r}, 
             \frac{n_1+n_2}{a}, ... \right] & \\
& \vdots & \\
I_{k+q+1}: & \left[ ..., \frac{n_1+(q-1)(n_1+n_2)}{(q-1)a+r}, \frac{n_1+q(n_1+n_2)}{b},
             \frac{n_1+n_2}{a}, ... \right] & \\
\end{eqnarray*}
Concluding that $a$ and $b$ are denominators of consecutive Farey fractions.
\end{proof}

Let $v_1 = (p_1,q_1,r_1)$ and $v_2 = (p_2,q_2,r_2)$ be the endpoints of an edge of a 
subtriangle in the $n^{th}$ stage of the Farey partitioning.  If $c = gcd(r_1,r_2)$, 
$r_1=cr'_1$, and $r_2=cr'_2$ then by the previous lemma we know there exist 
consecutive Farey fractions $n_1/r'_1$ and $n_2/r'_2$.  Looking at this edge at the
next iteration we have
\begin{eqnarray*}
& I_n     = \left[ \frac{1}{c}(\frac{p_1}{r'_1}, \frac{q_1}{r'_1}), 
                   \frac{1}{c}(\frac{p_2}{r'_2}, \frac{q_2}{r'_2}) \right] & \\
& I_{n+1} = \left[ \frac{1}{c}(\frac{p_1}{r'_1}, \frac{q_1}{r'_1}), 
                   \frac{1}{c}(\frac{p_1+p_2}{r'_1+r'_2}, \frac{q_1+q_2}{r'_1+r'_2}),
                   \frac{1}{c}(\frac{p_2}{r'_2}, \frac{q_2}{r'_2}) \right] & \\
\end{eqnarray*}
Now looking at what happens in the one-dimensional case
\begin{eqnarray*}
& I_k     = \left[ \frac{n_1}{r'_1}, \frac{n_2}{r'_2} \right] & \\
& I_{k+1} = \left[ \frac{n_1}{r'_1}, \frac{n_1 + n_2}{r'_1 + r'_2}, 
                   \frac{n_2}{r'_2} \right] & \\
\end{eqnarray*}
we see that the distances from endpoints to Farey sums are proportional by a factor of
$d(v_1,v_2)/d(n_1/r'_1,n_2/r'_2)$.  Since this holds for all combinations of Farey 
sums we must have
\[
\delta((1-x)v_1 + xv_2)  = 
     \frac{ ?((1-x)\frac{n_1}{r'_1} + x\frac{n_2}{r'_2}) - ?(\frac{n_1}{r'_1})}
          { ?(\frac{n_2}{r'_2}) - ?(\frac{n_1}{r'_1}) } \delta(v_2)
     +
     \frac{ ?(\frac{n_2}{r'_2}) - ?((1-x)\frac{n_1}{r'_1} + x\frac{n_2}{r'_2})}
          { ?(\frac{n_2}{r'_2}) - ?(\frac{n_1}{r'_1}) } \delta(v_1).
\]
Double checking with $v_1 = (0,0), v_2 = (1,0), n_1/r'_1 = 0/1, n_2 = r'_2 = 1/1$ we
find
\begin{eqnarray*}
\delta(x,0) &=& \delta( (1-x)(0,0) + x(1,0) ) \\
&=&
     \frac{ ?((1-x)\frac{0}{1} + x\frac{1}{1}) - ?(\frac{0}{1})}
          { ?(\frac{1}{1}) - ?(\frac{0}{1}) } \delta(1,0)
     +
     \frac{ ?(\frac{1}{1}) - ?((1-x)\frac{0}{1} + x\frac{1}{1}) }
          { ?(\frac{1}{1}) - ?(\frac{0}{1}) } \delta(0,0) \\
&=&
     \frac{ ?(x) - 0 }{1 - 0} (1,0) + \frac{ 1 - ?(x) }{ 1 - 0 } (0,0) \\
&=&  (?(x),0).
\end{eqnarray*}

\section{Farey Iteration as Multidimensional Continued Fraction}

We now want to show that $\delta$ sends a set of cubic irrational points to a set of 
rationals.  With this in mind, we first show that if $(\alpha,\beta) \in \Tri$ has 
an eventually periodic Farey sequence then $\alpha$ and $\beta$ are at worst cubic
irrational.  Later, we will show that if $(\alpha,\beta) \in \Tri$ has an 
eventually periodic Bary sequence then $\alpha$ and $\beta$ are rational.  Because
the Bary sequence of $\delta(\alpha,\beta)$ will equal the Farey sequence of 
$(\alpha,\beta)$ we can conclude that $\delta$ does send a natural class of cubic 
irrationals to a natural class of rationals.  To get there, we first define the 
revised Farey sequence.

We define the Farey sequence $\{ a_k(i_k) \}$ such that $a_k(i_k)$ denotes the 
action of performing $a_k$ type $i_k$ operations.  Where the six types of 
transformations on an arbitrary subtriangle $\langle v_1, v_2, v_3 \rangle$ are

\begin{figure}[hbtp]
\centering
\includegraphics[height=2.5in]{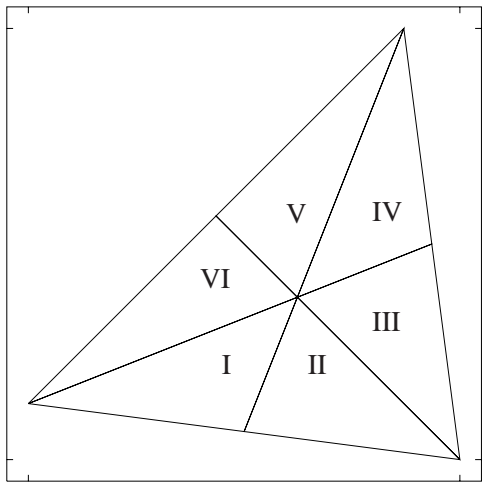}
\caption{The six Farey subtriangles.}
\end{figure}

\begin{enumerate}
\centering
\item $\langle v_1, 
        v_1 \fsum v_2, v_1 \fsum v_2 \fsum v_3 \rangle$
\item $\langle v_1, 
        v_1 \fsum v_3, v_1 \fsum v_2 \fsum v_3 \rangle$
\item $\langle v_2, 
        v_1 \fsum v_2, v_1 \fsum v_2 \fsum v_3 \rangle$
\item $\langle v_2, 
        v_2 \fsum v_3, v_1 \fsum v_2 \fsum v_3 \rangle$
\item $\langle v_3, 
        v_1 \fsum v_3, v_1 \fsum v_2 \fsum v_3 \rangle$
\item $\langle v_3, 
        v_2 \fsum v_3, v_1 \fsum v_2 \fsum v_3 \rangle$
\end{enumerate}
In proving singularness we are only interested in the subtriangle produced by taking 
$L = a_k$ type I moves, which is
\[
\langle v_1, Lv_1 \fsum v_2, \frac{L(L+1)}{2}v_1 \fsum Lv_2 \fsum v_3 \rangle.
\]

Now we will cast our six transformations in terms of matrix multiplication.  Recall 
that we write the initial triangle as
\[
M_0 =
( v_1(0) ~ v_2(0) ~ v_3(0) )
=
\left(
\begin{array}{ccc}
0 & 1 & 1 \\
0 & 0 & 1 \\
1 & 1 & 1 \\
\end{array}
\right).
\]
And we write subsequent triangles as
\[
M_n = 
( v_1(n) ~ v_2(n) ~ v_3(n) )
=
\left(
\begin{array}{ccc}
p_1(n) & p_2(n) & p_3(n) \\
q_1(n) & q_2(n) & q_3(n) \\
r_1(n) & r_2(n) & r_3(n) \\
\end{array}
\right).
\]
So, our six transformations have the form:
\begin{enumerate}
\centering
\item 
$M_n = M_{n-1} 
\left(
\begin{array}{ccc}
1 & 1 & 1 \\
0 & 1 & 1 \\
0 & 0 & 1 \\
\end{array}
\right)$

\item 
$M_n = M_{n-1} 
\left(
\begin{array}{ccc}
1 & 1 & 1 \\
0 & 0 & 1 \\
0 & 1 & 1 \\
\end{array}
\right)$

\item 
$M_n = M_{n-1} 
\left(
\begin{array}{ccc}
0 & 1 & 1 \\
1 & 1 & 1 \\
0 & 0 & 1 \\
\end{array}
\right)$

\item 
$M_n = M_{n-1} 
\left(
\begin{array}{ccc}
0 & 0 & 1 \\
1 & 1 & 1 \\
0 & 1 & 1 \\
\end{array}
\right)$

\item 
$M_n = M_{n-1} 
\left(
\begin{array}{ccc}
0 & 1 & 1 \\
0 & 0 & 1 \\
1 & 1 & 1 \\
\end{array}
\right)$

\item 
$M_n = M_{n-1} 
\left(
\begin{array}{ccc}
0 & 0 & 1 \\
0 & 1 & 1 \\
1 & 1 & 1 \\
\end{array}
\right)$
\end{enumerate}

\begin{thm} \label{determinant}
All possible $M_n$ have determinant $\pm 1$.
\end{thm}

\begin{proof}
Since $det(M_0) = 1$ and the determinants of the transformation matrices are also 
plus or minus one then any product of these matrices must have determinant plus 
or minus one.
\end{proof}

\begin{thm}
The entries in each column vector of $M_n$ share no common factor.
\end{thm}

\begin{proof}
This proof is identical to that of theorem $2.4.2$.
\end{proof}

This proves that our matrices are calculating the correct Farey sum.

\subsection{Farey Periodicity Implies Cubic Irrationality}

\begin{thm}
Suppose that $\{ a_k(i_k) \}$ is an eventually periodic Farey sequence that converges 
to the point $(\alpha,\beta)$.  Then $\alpha$ and $\beta$ are algebraic numbers with 
$deg(\alpha) \le 3$, $deg(\beta) \le 3$ and
\[
dim_{\Q}\Q[\alpha,\beta] \le 3.
\]
\end{thm}

\begin{proof}
The proof of theorem \ref{cubic} holds for any three by three matrices
hence it holds for any product of these six transformation matrices.
\end{proof}

\subsection{Not All Infinite Farey Sequences Converge}

\begin{exe}
Consider the Farey sequence consisting solely of type VI moves.  Having looked at the
first few subtriangles in the sequence, $v_2$ and $v_3$ appear to converge to the
point $(0.3820, 0.2361)$, while $v_1$ remains fixed at the origin.

\begin{figure}[hbtp]
\centering
\includegraphics[height=2.5in]{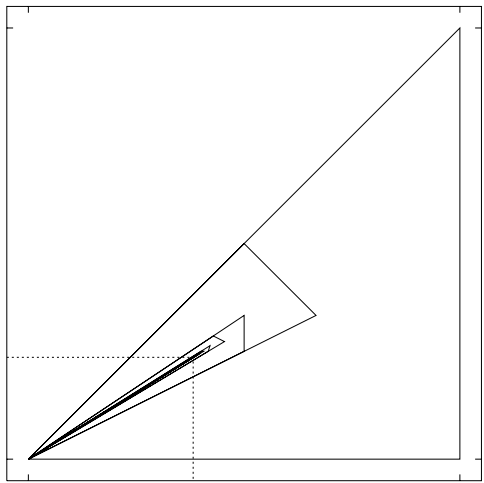}
\caption{An infinite Farey sequence that converges to a line segment.}
\end{figure}

We now prove this speculation.  Let $\{ f_n \}$ be the Fibonacci sequence, where 
$f_0 = 0, f_1 = 1$, and $f_n = f_{n-1}+f_{n-2}$.  Let 
$\langle v_1(n), v_2(n), v_3(n) \rangle$ denote the result of performing $n$ type VI 
transformations on $\Tri$.  We find the following pattern

\begin{center}
\begin{tabular}{rccc}
$0$:  & $v_1$  & $v_2$                         & $v_3$ \\
$1$:  & $v_1$  & $v_1 \fsum v_3$               & $v_1 \fsum v_2 \fsum v_3$ \\
$2$:  & $v_1$  & $2v_1 \fsum v_2 \fsum v_3$    & $3v_1 \fsum v_2 \fsum 2v_3$ \\
$3$:  & $v_1$  & $4v_1 \fsum v_2 \fsum 2v_3$   & $6v_1 \fsum 2v_2 \fsum 3v_3$ \\
$4$:  & $v_1$  & $7v_1 \fsum 2v_2 \fsum 3v_3$  & $11v_1 \fsum 3v_2 \fsum 5v_3$ \\
\vdots& \vdots & \vdots                        & \vdots \\
$n$:  & $v_1$  & $(f_{n+2}-1)v_1 \fsum f_{n-1}v_2 \fsum f_nv_3$
               & $(f_{n+3}-2)v_1 \fsum f_nv_2 \fsum f_{n+1}v_3$ \\
\end{tabular}
\end{center}

Taking the limit as $n$ goes to infinity gives the following
\begin{eqnarray*}
\mylim{n}{\infty} v_2(n) &=& \mylim{n}{\infty} \frac{f_{n+2}-1}{f_{n-1}}v_1 \fsum 
                             \frac{f_{n-1}}{f_{n-1}}v_2 \fsum \frac{f_n}{f_{n-1}}v_3\\
                         &=& \phi^3v_1 \fsum v_2 \fsum \phi v_3, 
\end{eqnarray*}
where $\phi = \mylim{n}{\infty} \frac{f_{n+1}}{f_n} = \frac{1+\sqrt{5}}{2}$.  
Similarly for $v_3$ we have
\begin{eqnarray*}
\mylim{n}{\infty} v_3(n) &=& \mylim{n}{\infty} \frac{f_{n+3}-2}{f_n}v_1 \fsum 
                             \frac{f_n}{f_n}v_2 \fsum \frac{f_{n+1}}{f_n}v_3 \\
                         &=& \phi^3v_1 \fsum v_2 \fsum \phi v_3.
\end{eqnarray*}
Now we use our empirical results to verify this limit.
\begin{eqnarray*}
\phi^3v_1 \fsum v_2 \fsum \phi v_3 
                 &=& (2+\sqrt{5})(0,0,1) ~ \fsum ~ (1,0,1) ~ \fsum ~ 
                     ((1+\sqrt{5})/2)(1,1,1)\\
                 &=& \left( \frac{3+\sqrt{5}}{2}, 
		            \frac{1+\sqrt{5}}{2},
			    \frac{7+3\sqrt{5}}{2} \right) \\
	         &=& \left( \frac{3+\sqrt{5}}{7+3\sqrt{5}},
			    \frac{1+\sqrt{5}}{7+3\sqrt{5}} \right)\\
		 &\approx& (0.3820, 0.2361).
\end{eqnarray*}

This is a surprising result.  We now know that there exist infinite Farey
sequences that converge to line segments.  From lemma \ref{sidelength} we know that
any infinite Bary sequence must converge to a point.  So, our continuous Farey-Bary 
map collapses entire line segments to single points in a continuous fashion.  For the
weighted Farey-Bary map there existed both infinite Farey and Bary sequences that
converged to line segments (infinite type I moves converged to the X-axis).
\end{exe}

\section{Iteration in the Barycentric Range}

To finish proving that $\delta$ sends a subset of cubic irrationals to a subset of 
rationals we must show that points with eventually periodic Bary sequences have 
rational coordinates.  And so we define the  Bary sequence $\{ \tilde{a}_k(i_k) \}$.  
The value of $\tilde{a}_k(i_k)$ denotes the action of applying $\tilde{a}_k$ type 
$i_k$ Bary transformations.  Recall that our six Bary transformations on a subtriangle
$\langle \tilde{v}_1, \tilde{v}_2, \tilde{v}_3 \rangle$ are
\begin{enumerate}
\centering
\item $\langle \tilde{v}_1, \tilde{v}_1 \bsum \tilde{v}_2, 
       \tilde{v}_1 \bsum \tilde{v}_2 \bsum \tilde{v}_3 \rangle$
\item $\langle \tilde{v}_1, \tilde{v}_1 \bsum \tilde{v}_3, 
       \tilde{v}_1 \bsum \tilde{v}_2 \bsum \tilde{v}_3 \rangle$
\item $\langle \tilde{v}_2, \tilde{v}_1 \bsum \tilde{v}_2, 
       \tilde{v}_1 \bsum \tilde{v}_2 \bsum \tilde{v}_3 \rangle$
\item $\langle \tilde{v}_2, \tilde{v}_2 \bsum \tilde{v}_3, 
       \tilde{v}_1 \bsum \tilde{v}_2 \bsum \tilde{v}_3 \rangle$
\item $\langle \tilde{v}_3, \tilde{v}_1 \bsum \tilde{v}_3, 
       \tilde{v}_1 \bsum \tilde{v}_2 \bsum \tilde{v}_3 \rangle$
\item $\langle \tilde{v}_3, \tilde{v}_2 \bsum \tilde{v}_3, 
       \tilde{v}_1 \bsum \tilde{v}_2 \bsum \tilde{v}_3 \rangle$
\end{enumerate}

Again we recast these operations in terms of matrices.

\begin{enumerate}
\centering
\item 
$M_n = M_{n-1} 
\left(
\begin{array}{ccc}
1 & 1/2 & 1/3 \\
0 & 1/2 & 1/3 \\
0 & 0 & 1/3 \\
\end{array}
\right)$

\item 
$M_n = M_{n-1} 
\left(
\begin{array}{ccc}
1 & 1/2 & 1/3 \\
0 & 0 & 1/3 \\
0 & 1/2 & 1/3 \\
\end{array}
\right)$

\item 
$M_n = M_{n-1} 
\left(
\begin{array}{ccc}
0 & 1/2 & 1/3 \\
1 & 1/2 & 1/3 \\
0 & 0 & 1/3 \\
\end{array}
\right)$

\item 
$M_n = M_{n-1} 
\left(
\begin{array}{ccc}
0 & 0 & 1/3 \\
1 & 1/2 & 1/3 \\
0 & 1/2 & 1/3 \\
\end{array}
\right)$

\item 
$M_n = M_{n-1} 
\left(
\begin{array}{ccc}
0 & 1/2 & 1/3 \\
0 & 0 & 1/3 \\
1 & 1/2 & 1/3 \\
\end{array}
\right)$

\item 
$M_n = M_{n-1} 
\left(
\begin{array}{ccc}
0 & 0 & 1/3 \\
0 & 1/2 & 1/3 \\
1 & 1/2 & 1/3 \\
\end{array}
\right)$
\end{enumerate}
Note that these three matrices are stochastic, meaning that their columns sum to
one.

\subsection{Bary Periodicity Implies Rationality}

\begin{thm}
If $\{ \tilde{a}_k(i_k) \}$ is an eventually periodic Bary sequence that 
converges to the point $(\alpha,\beta) \in \Tri$ then $\alpha$ and $\beta$ are 
rational.
\end{thm}

\begin{proof}
We refer back to the proof of theorem \ref{rational}, which applies here because all
the Bary transformation matrices are stochastic.
\end{proof}

\section{The Farey-Bary Analog of Singularness}

The key to Salem's proof that $?'(x) = 0$ a.e. is that if the limit
\[
\mylim{n}{\infty} \frac{ \text{length of interval in the range} }
                       { \text{length of interval in the domain} }
\]
exists and is finite then it vanishes almost everywhere.  We will prove a direct
generalization of this theorem for the continuous Farey-Bary map.  We will show
that if the following limit exists and is finite then at almost all points
\[
\mylim{n}{\infty} \frac{\text{area of subtriangle in the range}}
                       {\text{area of subtriangle in the domain}} = 0.
\]

\subsection{Areas of the Farey and Bary Subtriangles}

Recall our notation, that given a finite Farey sequence $\{ a_1(i_1),...,a_k(i_k) \}$
we defined
\[
\Tri_k = \{ (x,y) : \{a_1(i_1),...a_k(i_k)\} 
\text{ are the first $k$ terms in the Farey sequence} \}.
\]
\begin{figure}[hbtp]
\centering
\includegraphics[height=2.5in]{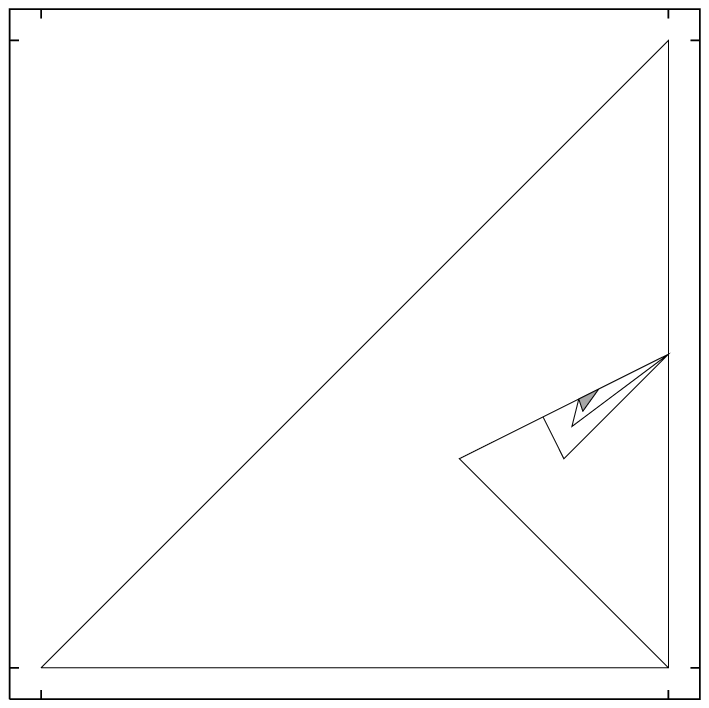}
\caption{$\Tri \{ 2(III), 1(I), 1(II) \}$.}
\end{figure}
The matrix representation $M_{s_k}$ of $\Tri_k$ is a product 
of the initial matrix $M_0$ and $s_k = a_1 + ... + a_k$ transformation matrices.

Using theorems \ref{area} and \ref{determinant} we have
\[
area( \Tri_k ) = \frac{1}{2} \cdot \frac{1}{r_1(k)r_2(k)r_3(k)}.
\]

At each stage of the Barycentric partitioning we are cutting the
area down by a factor of one sixth.  Hence
\[
area( \tilde{\Tri}_k ) = \frac{1}{2} \cdot \frac{1}{6^{s_k}}.
\]

Combining these two equations gives us the ratio
\[
\frac{area( \tilde{\Tri}_k )}
     {area( \Tri_k )} 
=
\frac{ r_1(k)r_2(k)r_3(k) }{6^{s_k}}.
\]
For the remainder of this paper, we will focus on showing that this ratio approaches 
zero almost everywhere.

\subsection{Almost everywhere $\limsup a_n(I) = \infty$}

In our quest to prove that $\delta$ is singular, we begin by showing that the number
of consecutive type I moves is unbounded almost everywhere.  

Given an arbitrary Farey sequence $\{ a_1(i_1), ... , a_k(i_k) \}$, let
$T = \Tri \{ a_1(i_1), ... , a_k(i_k) \}$ and
\[
M = (v_1 ~ v_2 ~ v_3) = 
\left(
\begin{array}{ccc}
p_1 & p_2 & p_3 \\
q_1 & q_2 & q_3 \\
r_1 & r_2 & r_3 
\end{array}
\right)
\]
be the matrix representation of $T$.  For any integer $L \ge 1$, define $T_L(I)$ to 
be the result of performing $L$ type I moves on $T$, thus
\[
T_L(I) = \langle v_1, Lv_1 \fsum v_2, 
         \frac{L(L+1)}{2}v_1 \fsum Lv_2 \fsum v_3 \rangle.
\]
Define $T_L = T - T_L(I)$.  $T_L$ is the collection of all points with 
$\{ a_1(i_1), ... , a_k(i_k) \}$ as the first $k$ terms of their Farey sequence and 
either $i_{k+1} \ne I$ or $a_{k+1} < L$.  

\begin{figure}[hbtp]
\centering
\includegraphics[width=\textwidth]{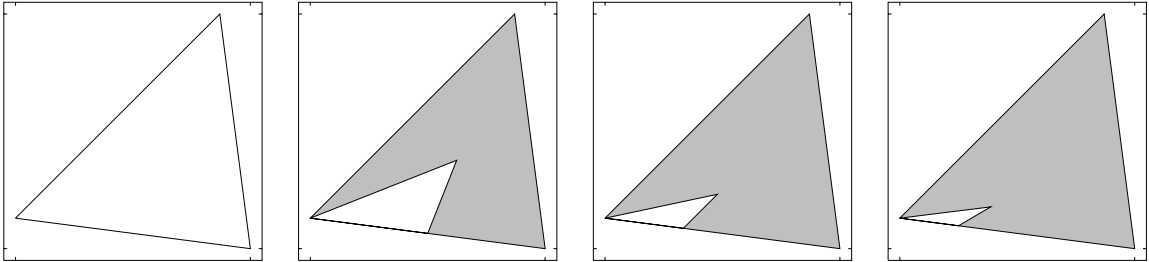}
\caption{$T_0,T_1,T_2,$ and $T_3$.}
\end{figure}

We now prove the following bound on $T_L$.

\begin{lem}
If $L \ge 1$ then 
\[
\frac{area(T_L)}{area(T)} \le \frac{8L^3-1}{8L^3}.
\]
\end{lem}

\begin{proof}
For ease of notation, set $x = r_1, y = r_2, z = r_3$.  Recall that
\[
2 \cdot area(T) = \frac{1}{xyz}.
\]
Now, we want to find a upper bound on the area of $T_L$.
\begin{eqnarray*}
2 \cdot area(T_L) &=&   \frac{1}{xyz} - \frac{1}{(x)(Lx+y)(\frac{L(L+1)}{2}x+Ly+z)}\\
                  &\le& \frac{1}{xyz} - \frac{1}{(x)(Ly+y)(\frac{L(L+1)}{2}z+Lz+z)}\\
                  &=&   \frac{1}{xyz} \left[ 1 - \frac{xyz}
                                                     {x(L+1)y(L/2+1)(L+1)z} \right]\\
                  &\le& \frac{1}{xyz} \left[ 1 - \frac{1}{(2L)(2L)(2L)} \right]\\
                  &=&   \frac{1}{xyz} \frac{8L^3 - 1}{8L^3}.
\end{eqnarray*}
\end{proof}

\begin{thm}
The set of $(\alpha,\beta) \in \Tri$ for which
\[
\mylimsup{n}{\infty} a_n(I) < \infty
\]
has measure zero.
\end{thm}

\begin{proof}
For each positive integer $N$, set
\[
M_N = \{ (\alpha,\beta) \in \Tri : 
          \forall n \ge 1, i_n \ne I \text{ or } a_n \le N \}.
\]
Since the union of all the $M_N$ is the set we want to show has measure zero if we can
show that $measure(M_N) = 0$ then we will be done.

Set
\[
M_N(1) = \{ (\alpha,\beta) \in \Tri : i_1 \ne I \text{ or } a_1 \le N\},
\]
and in general
\[
M_N(k) = \{ (\alpha,\beta) \in M_N(k-1) : i_k \ne I \text{ or } a_k \le N\}.
\]
Then we have a nested sequence of sets with
\[
M_N = \bigcap_{k=1}^\infty M_N(k).
\]

But this puts us into the language of the previous lemma.  Letting $L = N$, we have
\[
measure(M_N(k)) \le \frac{8N^3-1}{8N^3}measure(M_N(k-1)),
\]
and therefore
\[
measure(M_N) \le \prod_{k=1}^\infty \frac{8N^3-1}{8N^3} = 0.
\]
\end{proof}

\subsection{If the Limit Exists and is Finite then 
$\lim (area(\tilde{\Tri}_n) / area(\Tri_n)) = 0$}

We now prove the Farey-Bary analog of singularness.  If the following limit exists and
is finite then
\[
\mylim{n}{\infty} \frac{ \text{area of a subtriangle in the range} }
                       { \text{area of a subtriangle in the domain} } = 0
\]
almost everywhere.  This theorem is exactly what Salem used in his proof of 
singularness \cite{Salem, Derivative}.  In fact, we follow his proof closely.

\begin{thm}
Let $S = \{ (\alpha,\beta) \in \Tri : \mylimsup{n}{\infty} a_n(I) = \infty \}$.  For 
$(\alpha,\beta) \in S$, if
\[
\mylim{n}{\infty} \frac{area(\tilde{\Tri}\{a_1(i_1), ... , a_n(i_n)\})}
                       {area(\Tri\{a_1(i_1), ... ,a_n(i_n)\})}
\]
exists and is finite then it vanishes.
\end{thm}

\begin{proof}
Let $s_n = a_1 + ... + a_n$, and let
\begin{eqnarray*}
\rho_n &=& \frac{area(\tilde{\Tri}\{a_1(i_1),...,a_n(i_n)\})}
                     {area(\Tri\{a_1(i_1),...,a_n(i_n)\})} \\
       &=& \frac{r_1(n)r_2(n)r_3(n)}{6^{s_n}}.
\end{eqnarray*}
Now consider
\begin{eqnarray*}
\frac{\rho_n}{\rho_{n-1}} &=& \frac{r_1(n)r_2(n)r_3(n)}{6^{s_n}} 
                              \frac{6^{s_{n-1}}}{r_1(n-1)r_2(n-1)r_3(n-1)} \\
                          &=& \frac{r_1(n)r_2(n)r_3(n)}
			           {r_1(n-1)r_2(n-1)r_3(n-1)} \frac{1}{6^{a_n}}.
\end{eqnarray*}
For convenience set $x=r_1(n-1),y=r_2(n-1),$ and $z=r_3(n-1)$.  Assume that $i_n = I$ 
so that
\begin{eqnarray*}
\frac{\rho_n}{\rho_{n-1}} &=&   \frac{x(a_nx+y)(\frac{a_n(a_n+1)}{2}x+a_ny+z)}{xyz}
                                \frac{1}{6^{a_n}} \\
                          &\le& \frac{x(a_n+1)y(a_n+1)(a_n/2+1)z}{xyz}
                                \frac{1}{6^{a_n}}\\
                          &\le& \frac{8a_n^3}{6^{a_n}}.
\end{eqnarray*}

From the previous section, we know that $\mylimsup{n}{\infty} a_n(I) = \infty$, almost
everywhere.  Since the above denominator has a $6^{a_n}$ term while the numerator only
has a $a_n^3$ term, the entire ratio must approach zero
\[
\myliminf{n}{\infty} \frac{\rho_n}{\rho_{n-1}} = 0.
\]
If $\mylim{n}{\infty} \rho_n$ existed and were finite and different from $0$, then 
$\rho_n / \rho_{n-1}$ should tend necessarily to $1$.  Thus, if 
$\mylim{n}{\infty} \rho_n$ exists and is finite then $\mylim{n}{\infty} \rho_n = 0$.
\end{proof}

\subsection{Further Work}

In this paper, we defined two extensions of the Minkowski question-mark function.  Our
weighted map revealed a surprising asymmetry of the original Farey-Bary map, by 
setting the type I transformation apart from the other two transformations.  We 
defined the continuous Farey-Bary map and found reason to believe that it is a more 
natural extension of Minkowski's question-mark function.

There are a number of ways in which this work could be advanced.  Working directly 
from this paper, one could generalize the Farey-Bary maps to higher dimensions or 
place weights on the continuous Farey-Bary map.  A java applet that illustrates the 
weighted continuous Farey-Bary map can be found at 
http://wso.williams.edu/\~{}amarder/applets/ 

A more exciting problem is to find other two-dimensional Minkowski question-mark 
functions.  There are many multidimensional continued fractions.  For those that 
involve partitioning a triangle into any number of new subtriangles, can we define an 
analogous map to the Farey-Bary map?  If so, what are the properties of these new 
maps?

The more difficult and enlightening work would examine the link between the function
theoretic properties of the Farey-Bary map (singularness) and the number theoretic 
properties of the multidimensional continued fraction we defined using Farey 
iteration.

\nocite{*}                                             
\bibliography{biblio}                                  
\bibliographystyle{plain}                              

\end{document}